\documentclass[11pt]{article}

\usepackage{amsmath}

\usepackage{amssymb}
\usepackage{latexsym}
\usepackage{dsfont}
\begin{document}

\newtheorem{theo}{Theorem}
\newtheorem{prop}{Proposition}
\newtheorem{lem}{Lemma}
\newtheorem{cor}{Corollary}
\newtheorem{defi}{Definition}
\newtheorem{exple}{Example}
\newtheorem{Rk}{Remark}

\begin{center}
  {\Large \bf \sc Large and moderate deviations principles for 
\vspace{0.25cm}
\\
kernel estimators of the multivariate regression 
}\\
\vspace{1.cm}
Abdelkader Mokkadem \hspace{1.2cm} Mariane Pelletier \hspace{1.2cm} Baba Thiam

\vspace{0.3cm}
(mokkadem, pelletier, thiam)@math.uvsq.fr

\vspace{0.5cm}

{\small Universit\'e de Versailles-Saint-Quentin\\
D\'epartement de Math\'ematiques\\
45, Avenue des Etats-Unis\\
78035 Versailles Cedex\\
France }\\
\vspace{1.cm}

{\bf Abstract :}\ \parbox[t]{9.cm}{\small
In this paper, we prove large deviations principle for the Nadaraya-Watson estimator and for the semi-recursive kernel estimator of the regression in the multidimensional case. Under suitable conditions, we show that the rate function is a good rate function. We thus generalize the results already obtained in the unidimensional case for the Nadaraya-Watson estimator. Moreover, we give a moderate deviations principle for these two estimators. It turns out that the rate function obtained in the moderate deviations principle for the semi-recursive estimator is larger than the one obtained for the Nadaraya-Watson estimator.}
\vspace{2.cm}

\rm

\end{center}
\noindent
{\bf AMS Subj. Classification:} 62G08, 60F10\\
\noindent
{\bf Key words and phrases~:}\\
Nadaraya-Watson estimator~; Recursive kernel estimator~; Large deviations principle~; Moderate deviations principle 
\vspace{0.5cm}

\newpage
\section{Introduction}
Let $(X,Y)$, $(X_1,Y_1), \ldots,(X_n,Y_n)$ be a sequence of independent and identically distributed $\mathbb{R}^d\times \mathbb{R}^q$-valued random variables with probability density $f(x,y)$ with $\mathbb{E}|Y|<\infty$. Moreover, let $g(x)$ be the marginal density of $X$ and $r(x)=\mathbb{E}\left(Y\vert X=x\right)=m(x)/g(x)$ the regression of $Y$ on $X$. The purpose of this paper is to establish large and moderate deviations principles for the Nadaraya-Watson estimator and for the semi-recursive kernel estimator of the regression.\\

Let us first recall the concept of large and moderate deviations. A speed is a sequence $(\nu_n)$ of positive numbers going to infinity. A good rate function on $\mathbb{R}^m$ is a lower semicontinuous function $I:\mathbb{R}^m \rightarrow [0,\infty]$ such that, for each $\alpha<\infty$, the level set $\{x\in \mathbb{R}^m,\ \ I(x)\leq \alpha\}$ is a compact set. If the level sets of $I$ are only closed, then $I$ is said to be a rate function. A sequence $(Z_n)_{n \geq 1}$ of $\mathbb{R}^m$-valued random variables is said to satisfy a large deviations principle (LDP) with speed $(\nu_n)$ and rate function $I$ if:
\begin{eqnarray*}
\liminf_{n \to \infty}\nu_n^{-1}\log\mathbb{P}\left[Z_n\in U\right] & \geq &-\inf_{x \in U}I(x) \ \ \mbox{for every open subset U of}\ \ \mathbb{R}^m,\\
\limsup_{n \to  \infty}\nu_n^{-1}\log\mathbb{P}\left[Z_n\in V\right]  &\leq & 
-\inf_{x \in V}I(x)\ \ \mbox{for every closed subset V of}\ \ \mathbb{R}^m.
\end{eqnarray*}
Moreover, let $(v_n)$ be a nonrandom sequence that goes to infinity; if $(v_nZ_n)$ satisfies a LDP, then $(Z_n)$ is said to satisfy a moderate deviations principle (MDP).\\

The Nadaraya-Watson estimator (\cite{nadarayar}, \cite{watsonr}) of the regression function $r(x)$ is defined by
\begin{eqnarray}
\label{estimateregression}
r_n(x)=
\left\{
\begin{array}{ll}
\dfrac{m_n(x)}{g_n(x)} \ \ \mbox{if} \ \ g_n(x) \neq 0\\
0 \ \ \mbox{otherwise,} 
\end{array}
\right.
\end{eqnarray}
with
\begin{eqnarray*}
m_n(x) = \frac{1}{nh_n^d}\sum_{i=1}^nY_iK\left(\frac{x-X_i}{h_n}\right) \ \ \mbox{and} \ \ 
g_n(x)=\frac{1}{nh_n^d}\sum_{i=1}^nK\left(\frac{x-X_i}{h_n}\right),
\end{eqnarray*}
where the bandwidth $(h_n)$ is a positive sequence such that 
\begin{eqnarray}
\label{bandwidth}
\lim_{n \to \infty}h_n=0 \ \ \mbox{and} \ \  \lim_{n \to \infty}nh_n^d=\infty,
\end{eqnarray}
and the kernel $K$ a continuous function such that $\lim_{\|x\| \to\infty}K(x)=0$ and $\int_{\mathbb{R}^d}K(x)dx=1$. The weak and strong consistency of $r_n$ has been widely discussed by many authors; let us cite, among many others, Collomb \cite{colr}, Collomb and H\"ardle \cite{colharr}, Devroye \cite{devr}, Mack and Silverman \cite{macsilr} and Senoussi \cite{senr}. For other works on the consistency of $r_n$, the reader is refered to the monographs of Bosq \cite{bosr} and Prakasa Rao \cite{prar}. The large deviations behaviour of $r_n$ has been studied at first by Louani \cite{lour}, and then by Joutard \cite{jour} in the univariate framework. Moderate deviations principles have been obtained by Worms \cite{wormr} in the particular case $Y=r(X)+\varepsilon$ with $\varepsilon$ and $X$ independent. The first aim of this paper is to generalize these large and moderate deviations results.\\

The approach used by Louani \cite{lour} and Joutard \cite{jour} to study the large deviations behaviour of $r_n$ is to note that, if $d=q=1$ and if the kernel is positive, then, for all $\delta>0$, 
\begin{eqnarray*}
\mathbb{P}\left[r_n(x)-r(x)\geq \delta\right] & = & \mathbb{P}\left[\frac{1}{nh_n}\sum_{j=1}^n\left[Y_j-r(x)-\delta\right]K\left(\frac{x-X_j}{h_n}\right) \geq 0\right].
\end{eqnarray*}
Obviously, their approach can not be extended to the multivariate framework. Thus, to study the large deviations behaviour of $r_n$, our approach is totally different. We first establish a large deviations principle for the sequence $\left(m_n(x),g_n(x)\right)$, and then show how the large deviations behaviour of $r_n$ can be deduced. More precisely, for $x \in \mathbb{R}^d$, let $\Psi_x$ be the function defined for any $(u,v)\in \mathbb{R}^q\times\mathbb{R}$ by
\begin{eqnarray*}
\Psi_x(u,v) & = & \int_{\mathbb{R}^d\times \mathbb{R}^q}\left(e^{(\langle u,y\rangle+v)K(z)}-1\right)f(x,y)dzdy,
\end{eqnarray*}
(where $\langle u,y\rangle$ denotes the scalar product of $u$ and $y$) and let $I_x$ be the Fenchel-Legendre transform of $\Psi_x$. We give conditions ensuring that the sequence $\left(r_n(x)\right)$ satisfies a LDP with speed $(nh_n^d)$ and good rate function $J$ defined, for any $s\in \mathbb{R}^q$, by 
\begin{eqnarray*}
J(s)=\inf_{t \in \mathbb{R}}I_x(st,t).
\end{eqnarray*} 

Concerning the moderate deviations behaviour of the Nadaraya-Watson estimator, we prove that, for any positive sequence $(v_n)$ such that
\begin{eqnarray}
\label{condition}
\lim_{n \to \infty}v_n=\infty , \ \ \lim_{n \to \infty}\frac{v_n^2}{nh_n^d}=0,\ \ \mbox{and}\ \ \lim_{n \to \infty}v_nh_n^p=0,
\end{eqnarray} 
(where $p$ denotes the order of the kernel $K$) the sequence $\left(v_n\left[r_n(x)-r(x)\right]\right)$ satisfies a LDP with speed $\left(nh_n^d/v_n^2\right)$ and good rate function $G_x$ defined for all $v \in \mathbb{R}^q$ by
\begin{eqnarray}
\label{rateregression}
G_x(v) & = & \frac{g(x)}{2\int_{\mathbb{R}^d}K^2(z)dz}v^T \Sigma_x^{-1} v,
\end{eqnarray}
where $\Sigma_x$ denotes the $q\times q$ covariance matrix $V\left(Y\vert X=x\right)$. 
Let us note that, in the case the model $Y=r(X)+\varepsilon$ (with $X$ and $\varepsilon$ independent) is considered, the matrix $\Sigma_x$ is the covariance matrix of $\varepsilon$ and does depend on $x$; we then find the MDP proved in Worms \cite{wormr} again.\\

A semi-recursive version of the Nadaraya-Watson estimator \eqref{estimateregression} is defined as 
\begin{eqnarray}
\label{estimatebisregression}
\tilde r_n(x)=
\left\{
\begin{array}{ll}
\dfrac{\tilde m_n(x)}{\tilde g_n(x)} \ \ \mbox{if} \ \ \tilde g_n(x) \neq 0\\
0 \ \ \mbox{otherwise},
\end{array}
\right.
\end{eqnarray}
where 
\begin{eqnarray*}
\tilde m_n(x)=\frac{1}{n}\sum_{i=1}^n\frac{Y_i}{h_i^d}K\left(\frac{x-X_i}{h_i}\right) \ \ \mbox{and}\ \ \tilde g_n(x)=\frac{1}{n}\sum_{i=1}^n\frac{1}{h_i^d}K\left(\frac{x-X_i}{h_i}\right).
\end{eqnarray*}
Weak conditions for various forms of consistency of $\tilde r_n$ have been obtained by Ahmad and Lin \cite{ahmadlinr} and Devroye and Wagner \cite{devroyewagnerr}. Roussas \cite{roussasr} studied its almost sure convergence rate. The second aim of this paper is to establish the large and moderate deviations behaviour of $\tilde r_n$.\\

It turns out that the rate function that appears in the LDP is much more complex to explicit in the case the semi-recursive kernel regression estimator is considered than in the case the Nadaraya-Watson estimator is used. That is the reason why we only consider bandwidths defined as $(h_n)=(cn^{-a})$ with $c>0$ and $a \in ]0,1/d[$ (instead of bandwidths satisfying \eqref{bandwidth}). For $x \in \mathbb{R}^d$, let $\tilde\Psi_{a,x}$ be the function defined for all $(u,v)\in \mathbb{R}^q\times\mathbb{R}$ by
\begin{eqnarray*}
\tilde \Psi_{a,x}(u,v) & = & \int_{[0,1]\times\mathbb{R}^d\times \mathbb{R}^q}s^{-ad}\left(e^{s^{ad}(\langle u,y\rangle+v)K(z)}-1\right)f(x,y)dsdzdy,
\end{eqnarray*}
and let $\tilde I_{a,x}$ be the Fenchel-Legendre transform of $\tilde \Psi_{a,x}$. We give conditions ensuring that the sequence $\left(\tilde r_n(x)\right)$ satisfies a LDP with speed $(nh_n^d)$ and good rate function $\tilde J_a$ defined, for any $s\in \mathbb{R}^q$, by
\begin{eqnarray*}
\tilde J_a(s)=\inf_{t \in \mathbb{R}}\tilde I_{a,x}(st,t).
\end{eqnarray*}

To establish the moderate deviations behaviour of $\tilde r_n$, we consider bandwidths $(h_n)$ which vary regularly with exponent $(-a)$, $a \in ]0,1/d[$. We prove that, for any positive sequence $(v_n)$ satisfying \eqref{condition}, the sequence $\left(v_n\left[\tilde r_n(x)-r(x)\right]\right)$ satisfies a LDP with speed $\left(nh_n^d/v_n^2\right)$ and good rate function defined for all $v \in \mathbb{R}^q$ by
\begin{eqnarray}
\label{rateprimregression}
\tilde G_{a,x}(v) & = & \frac{(1+ad)g(x)}{2\int_{\mathbb{R}^d}K^2(z)dz}v^T \Sigma_x^{-1}v.
\end{eqnarray}
Let us underline that, because of the factor $(1+ad)$ which is present in \eqref{rateprimregression} but not in \eqref{rateregression}, the rate function obtained in the MDP in the case the semi-recursive estimator is used is larger than the one which appears in the case the Nadaraya-Watson kernel estimator is considered; this means that the semi-recursive estimator $\tilde r_n(x)$ is more concentrated around $r(x)$ than the Nadaraya-Watson estimator.\\

Our main results are stated in Section $2$, whereas Section $3$ is devoted to the proofs.
\section{Assumptions and Main Results}
We shall use the following notations.
\begin{itemize}
\item $\mathcal{D}(\mathcal{F})=\left\{x, \ \mathcal{F}(x)<\infty\right\}$ denotes the domain of a function $\mathcal{F}$ and $\overset{\circ}{\mathcal{D}}(\mathcal{F})$ is the interior domain of $\mathcal{F}$.
\item $\|x\|$ is the euclidean norm of $x$.
\item $\lambda$ is the Lebesgue measure.
\item $a\wedge b=\min\{a,b\}$. 
\item $\vec{0}=\left(0,\ldots,0\right) \in \mathbb{R}^q$.
\end{itemize}

The large and moderate deviations behaviours of the Nadaraya-Watson estimator $r_n$ are given in Section $2.1$, whereas the ones of the semi-recursive kernel estimator $\tilde r_n$ are stated in Section $2.2$.
\subsection{Large and moderate deviations principles for the Nadaraya-Watson estimator}
The assumptions required for the LDP of the Nadaraya-Watson estimator are the following. 
\begin{description}
\item (A1) $K:\mathbb{R}^d\rightarrow \mathbb{R}$ is a bounded and integrable function, $\int_{\mathbb{R}^d}K(z)dz=1$ and $\lim_{\|z\| \to \infty}K(z)=0$.
\item (A2) For any $u \in \mathbb{R}^q$, $t\mapsto\displaystyle \int_{\mathbb{R}^q}e^{\langle u,y\rangle }f(t,y)dy$ is continuous at $x$ and bounded.    
\end{description}
\paragraph{Comments}
\begin{itemize}
\item Notice that (A2) implies that the density $g$ is continuous at $x$ and bounded. 
\item In the model $Y=r(X)+\varepsilon$ with $\varepsilon$ and $X$ independent, let $h$ be the probability density of $\varepsilon$. Then
\begin{eqnarray*}
f(t,y)& = & g(t)h\left(y-r(t)\right)\\
\int_{\mathbb{R}^q}\|y\|f(t,y)dy &= & g(t)\int_{\mathbb{R}^q}\|y+r(t)\|h(y)dy \\ 
\int_{\mathbb{R}^q}e^{\langle u,y\rangle }f(t,y)dy & = & g(t)e^{\langle u,r(t)\rangle }\int_{\mathbb{R}^q}e^{\langle u,y\rangle }h(y)dy.
\end{eqnarray*}
Thus, (A2) can be translated as assumptions on $g$ and $r$ and on the moments of $\varepsilon$.
\item As it can be seen from the proofs, the boundness assumption in (A2) is useless if $K$ has a compact support.
\item The boundness of the function $t\mapsto\displaystyle\int_{\mathbb{R}^q}e^{\langle u,y\rangle }f(t,y)dy$ for any $u \in \mathbb{R}^q$ implies that 
\begin{eqnarray}
\label{borne}
\forall m\geq 0, \forall \rho\geq 0 \ \ \mbox{the function}\ \ t\mapsto\int_{\mathbb{R}^q}\|y\|^me^{\rho\|y\|}f(t,y)dy \ \ \mbox{ is bounded.}
\end{eqnarray}

\paragraph{Proof}
It suffices to prove that the function $t\mapsto\displaystyle\int_{\mathbb{R}^q}e^{\rho \|y\|}f(t,y)dy$ is bounded for any $\rho>0$. Set $y=\left(y_1,\ldots,y_q\right)$, we first note that 
\begin{eqnarray*}
\int_{\mathbb{R}^q}e^{q\rho|y_j|}f(t,y)dy & \leq & 
\int_{\left\{y_j\geq 0\right\}}e^{q\rho y_j}f(t,y)dy+\int_{\left\{y_j < 0\right\}}e^{-q\rho y_j}f(t,y)dy\\ & \leq & 
\int_{\mathbb{R}^q}e^{q\rho y_j}f(t,y)dy+\int_{\mathbb{R}^q}e^{-q\rho y_j}f(t,y)dy.
\end{eqnarray*}
Now, we have
\begin{eqnarray*}
\lefteqn{\int_{\mathbb{R}^q}e^{\rho \|y\|}f(t,y)dy}\\ & \leq & 
\int_{\mathbb{R}^q}e^{\rho |y_1|+\cdots+\rho|y_q|}f(t,y)dy\\& \leq & 
\left(\int_{\mathbb{R}^q}e^{q\rho|y_1|}f(t,y)dy\ldots\int_{\mathbb{R}^q}e^{q\rho|y_q|}f(t,y)dy\right)^{\frac{1}{q}} \ \ \mbox{by the generalized H\"older inequality.}\\
& \leq &
\left(\left(\int_{\mathbb{R}^q}e^{q\rho y_1}f(t,y)dy+\int_{\mathbb{R}^q}e^{-q\rho y_1}f(t,y)dy\right)\ldots\left(\int_{\mathbb{R}^q}e^{q\rho y_q}f(t,y)dy+\int_{\mathbb{R}^q}e^{-q\rho y_q}f(t,y)dy\right)\right)^{\frac{1}{q}}
\end{eqnarray*}
which is bounded.

\end{itemize}

Before stating our results, we need to introduce the rate function for the LDP of the Nadaraya-Watson estimator. Let $\Psi_x:\mathbb{R}^q\times\mathbb{R}\rightarrow\mathbb{R}$ and $I_x,\hat{I}_x:\mathbb{R}^q\times\mathbb{R}\rightarrow\mathbb{R}$ be the functions defined as follows: 
\begin{eqnarray}
\label{psiregression}
 \Psi_x(u,v) & = & \int_{\mathbb{R}^d\times \mathbb{R}^q}\left(e^{(\langle u,y\rangle+v)K(z)}-1\right)f(x,y)dzdy,\\
\label{tauxregression}
I_x(t_1,t_2)  &= & \sup_{(u,v) \in \mathbb{R}^q\times\mathbb{R}}\left\{\langle u,t_1\rangle+vt_2-\Psi_x(u,v)\right\},\\
\hat{I}_x(s,t) & = & I_x(st,t).
\label{tauxnadwat}
\end{eqnarray}
Moreover, for any $s \in \mathbb{R}^q$, set  
\begin{eqnarray*}
J^*(s) & = & \inf_{t \in \mathbb{R}^*}I_x(st,t) \\
  & = & \inf_{t \in \mathbb{R}^*}\hat I_x(s,t),\\
J(s) & = & J^*(s)\wedge I_x(\vec{0},0)\\
        & = & \inf_{t \in \mathbb{R}}\hat I_x(s,t).
\end{eqnarray*} 
To prove that $J$ is a rate function, we need to assume that the following condition (C) is fulfilled.
\begin{description}
\item (C) $\inf_{s \in \mathbb{R}^q}I_x(s,0)=I_x(\vec{0},0)$.
\end{description}
Before stating the properties of the function $J$, let us give some cases when Condition (C) is satisfied (under Assumptions (A1) and (A2)).
\paragraph{Example 1:} \textbf{Nonnegative kernel}\\
Condition (C) is satisfied when $K$ is nonnegative since, in this case, $I_x(s,0)=+\infty$ for any $s\neq \vec{0}$, (this is stated in Proposition \ref{convex1} of Section $3$).
\paragraph{Example 2:} \textbf{Model with symmetry}\\
Condition (C) holds when $f$ is symmetric in each coordinate of the second variable $y\in \mathbb{R}^q$. As a matter of fact, for a diagonal $q\times q$ matrix $A$ such that $A_{ii}=\pm1$, observe that
\begin{eqnarray*}
\Psi_x(Au,v) & = & \int_{\mathbb{R}^d\times \mathbb{R}^q}\left(e^{(\langle Au,y\rangle+v)K(z)}-1\right)f(x,y)dzdy\\  & = & 
 \int_{\mathbb{R}^d\times \mathbb{R}^q}\left(e^{(\langle u,Ay\rangle+v)K(z)}-1\right)f(x,y)dzdy\\
& = & 
 \int_{\mathbb{R}^d\times \mathbb{R}^q}\left(e^{(\langle u,y'\rangle+v)K(z)}-1\right)f(x,A^{-1}y')dzdy' \\& = &\Psi_x(u,v).
\end{eqnarray*}
For any given $s \in \mathbb{R}^q$, set 
\begin{eqnarray*}
\mathcal{U}_s=\left\{u \in \mathbb{R}^q, \ \langle u,s\rangle \geq 0\right\}.
\end{eqnarray*}
We have,
\begin{eqnarray*}
\sup_{u,v}\left(-\Psi_x(u,v)\right) & = & \sup_{u \in \mathcal{U}_s, \ v\in \mathbb{R}}\left(-\Psi_x(u,v)\right).
\end{eqnarray*}
Now, for any $u \in \mathcal{U}_s$ and $v  \in \mathbb{R}$,
\begin{eqnarray*}
\langle u,s\rangle-\Psi_x(u,v)\geq -\Psi_x(u,v),
\end{eqnarray*}
so that
\begin{eqnarray*}
\sup_{u,v}\left\{\langle u,s\rangle-\Psi_x(u,v)\right\} \geq
\sup_{u \in \mathcal{U}_s,\ v\in \mathbb{R}}\left(-\Psi_x(u,v)\right),
\end{eqnarray*}
and thus,
\begin{eqnarray*}
I_x(\vec{0},0) \leq I_x(s,0) \ \ \forall s\in \mathbb{R}^q,
\end{eqnarray*}
so that Condition (C) follows.
\paragraph{Example 3:}\textbf{A negative kernel without symmetry assumption on \boldmath$f$, and for $d=q=1$}\\
 If the kernel $K$ can be written as $\displaystyle K=\mathds{1}_{D}-\mathds{1}_{D'}$ where $D$ and $D'$ are two subsets of $\mathbb{R}$ such that $D\cap D' =\emptyset$ and $\lambda(D)-\lambda(D')=1$, then Condition (C) holds. As a matter of fact, we then have 
\begin{eqnarray*}
 \Psi_x(u,v) & = & \int_{\mathbb{R}\times \mathbb{R}}\left(e^{( uy+v)K(z)}-1\right)f(x,y)dzdy \\ & = & 
\int_{ D\times\mathbb{R}}\left(e^{ uy+v}-1\right)f(x,y)dzdy    +\int_{D'\times\mathbb{R} }\left(e^{ -uy-v}-1\right)f(x,y)dzdy \\ & = & 
e^v\lambda(D)\int_{\mathbb{R}}e^{uy}f(x,y)dy -\left[\lambda(D)+\lambda(D')\right]g(x)  +e^{-v}\lambda(D')\int_{\mathbb{R}}e^{-uy}f(x,y)dy.
\end{eqnarray*}
Now, let $M_x$ denote  the Laplace transform of $f(x,\cdot)$, then 
\begin{eqnarray*}
 \Psi_x(u,v) & = & e^v\lambda(D)M_x(u)  +e^{-v}\lambda(D')M_x(-u)-\left[\lambda(D)+\lambda(D')\right]g(x).
\end{eqnarray*}
For any given $u$, it can easily be seen that the infimum of $\Psi_x(u,\cdot)$ is reached at 
\begin{eqnarray*}
v_0=\log\sqrt{\frac{\lambda(D')M_x(-u)}{\lambda(D)M_x(u)}},
\end{eqnarray*}

and
\begin{eqnarray*}
 \Psi_x(u,v_0) & = &2\sqrt{\lambda(D)\lambda(D')}\sqrt{M_x(u)M_x(-u)}-\left[\lambda(D)+\lambda(D')\right]g(x).
\end{eqnarray*}
Observe that 
\begin{eqnarray*}
\Psi_x(u,v_0)=\Psi_x(-u,v_0),
\end{eqnarray*}
and thus
\begin{eqnarray*}
\sup_{u}\left(-\Psi_x(u,v_0)\right) & = & \sup_{u\geq 0}\left(-\Psi_x(u,v_0)\right)\\ & = & 
\sup_{u\leq 0}\left(-\Psi_x(u,v_0)\right)\\ & = & 
I_x(0,0).
\end{eqnarray*}
Now, if $s\geq 0$, we have for any $u\geq 0$ 
\begin{eqnarray*}
us- \Psi_x(u,v_0)\geq  -\Psi_x(u,v_0), 
\end{eqnarray*}
and thus
\begin{eqnarray*}
  I_x(s,0)& \geq & I_x(0,0)\ \forall s\geq 0.
\end{eqnarray*}
Proceeding in the same way for $s < 0$, we obtain Condition (C).\\
Such an example of a four order kernel is 
$\displaystyle K=\mathds{1}_{[-a,a]}-\mathds{1}_{[-b,-a[\cup ]a,b]}$, with
\begin{eqnarray*} 
 a & = & \frac{1}{6}\sqrt[3]{2}+\frac{1}{12}\left(\sqrt[3]{2}\right)^2+\frac{1}{3}\\
b & = & \frac{1}{3}\sqrt[3]{2}+\frac{1}{6}\left(\sqrt[3]{2}\right)^2+\frac{1}{6}.
\end{eqnarray*}

Let us now give the properties of the function $J$.
\begin{prop}\label{barca}
Assume that (A1), (A2) and (C) hold. Then,
\begin{description}
\item (i) $J$ is a rate function on $\mathbb{R}^q$. More precisely, for $\alpha\in \mathbb{R}$,
\begin{itemize}
\item if $\alpha<I_x(\vec{0},0)$, then $\left\{J(s)\leq \alpha\right \}$ is compact.
\item if $\alpha\geq I_x(\vec{0},0)$, then $\left\{J(s)\leq \alpha\right \}=\mathbb{R}^q$.
\end{itemize}
\item (ii) If $I_x(\vec{0},0)=\infty$, then $J$ is a good rate function on $\mathbb{R}^q$ and $J=J^*$.
\item (iii) If $J^*(s)<\infty$, then $J(s)=J^*(s)$.
\item (iv) If $\alpha<I_x(\vec{0},0)$, then $\left\{J^*(s)\leq \alpha\right \}=\left\{ J(s)\leq \alpha\right \}$. 

\end{description}

\end{prop}

\begin{Rk}\label{jii}
In view of the definition of $J$ and $J^*$, and of Proposition \ref{barca} (iii), we have:
\begin{eqnarray*}
J(s) & = & 
\left\{
\begin{array}{ll}
J^*(s)\ \ \mbox{if} \ \ J^*(s)<\infty\\
I_x(\vec{0},0)\ \ \mbox{if} \ \ J^*(s)=\infty.
\end{array}
\right.
\end{eqnarray*}
\end{Rk}
Let us now state the LDP for the Nadaraya-Watson estimator.

 \begin{theo}\label{ldpregression}\textbf{(Pointwise LDP for the Nadaraya-Watson estimator)}$ $ \\
Assume that (A1), (A2) and (C) hold, and that $(h_n)$ satisfies the conditions in \eqref{bandwidth}. Then, for any open subset $U$ of $\mathbb{R}^q$,  
\begin{eqnarray*}
 \liminf_{n \to \infty}\frac{1}{nh_n^d} \log\mathbb{P}\left[r_n(x)\in U\right] \geq -\inf_{s\in U}J^*(s), 
\end{eqnarray*}
and for any closed subset $V$ of $\mathbb{R}^q$,
\begin{eqnarray*}
 \limsup_{n \to \infty}\frac{1}{nh_n^d}\log \mathbb{P}\left[r_n(x)\in V\right] \leq -\inf_{s\in V}J(s).
\end{eqnarray*}
\end{theo}
\paragraph{Comments.}
\begin{description}
\item 1) Set $E=\left\{J^*(s)<\infty\right\}$. For any open subset $U$ of $\mathbb{R}^q$ such that $U\cap E\neq\emptyset$, we have
\begin{eqnarray*}
 \liminf_{n \to \infty}\frac{1}{nh_n^d} \log\mathbb{P}\left[r_n(x)\in U\right] \geq -\inf_{s\in U}J(s).
\end{eqnarray*}
\item 2) If $I_x$ is finite in a neighbourhood of $(\vec{0},0)$, then $J^*$ is finite everywhere and by Proposition \ref{barca} (iii), $\displaystyle J(s)=J^*(s)<\infty \ \ \forall s$ ; thus $(r_n)$ satisfies a LDP with speed $(nh_n^d)$ and rate function $J$.\\
Of course, this does not hold for nonnegative kernel since in this case $I_x(s,0)=+\infty$ for any $s$ (see Proposition \ref{convex1} in Section 3). However, it can hold for kernels which take negative values. For example, consider the previous Example $3$, and assume $f(x,y)$ is symmetric in $y$ ; in this case $\displaystyle M_x(u)=M_x(-u)$. The equation
\begin{eqnarray*}
\frac{\partial\Psi_x}{\partial v}(u,v)  =  \left[\lambda(D)e^{v}  -\lambda(D')e^{-v}\right]M_x(u)   = 0
\end{eqnarray*}
has solution $v_0= \log\sqrt{\dfrac{\lambda(D')}{\lambda(D)}}$ independent from $u$. Moreover, $M'$ is continuous and has range $\mathbb{R}$, thus, there exists $u_0$ such that $M'(u_0)=0$. This implies that the equation
\begin{eqnarray*}
\frac{\partial\Psi_x}{\partial u} (u,v) =  \left[\lambda(D)e^{v}  +\lambda(D')e^{-v}\right]M'_x(u)   = 0
\end{eqnarray*}
has a solution $u_0$ independent from $v$. Thus $(0,0)$ is in the range of $\nabla\Psi_x$. It follows from Proposition \ref{convex1} Section $3$ that $I_x$ is finite in a neighbourhood of $(0,0)$.
\item 3) When $I_x(\vec{0},0)=\infty$, it follows from Proposition \ref{barca} and Theorem \ref{ldpregression} that $(r_n)$ satisfies a LDP with speed $(nh_n^d)$ and good rate function $J$. 
\end{description}

In the case $K$ is a nonnegative kernel whose support has an infinity measure, we will show in Proposition \ref{convex1} that $I_x(\vec{0},0)=\infty$. We have thus the following corollary.
\begin{cor}\label{nonnegative}
Let the assumptions of Theorem \ref{ldpregression} hold. If $K$ is a nonnegative kernel such that $\lambda\left(\left\{x\in\mathbb{R}^d, K(x)>0\right\}\right)=\infty$, then the sequence $(r_n)$ satisfies a LDP with speed $(nh_n^d)$ and good rate function $J$.
\end{cor}

This corollary is an extension of the results of Louani \cite{lour} and Joutard \cite{jour} to the multivariate framework (and to the case the kernel $K$ may vanish). Moreover, it proves that the rate function that appears in their large deviations results is in fact a good rate function.\\

To establish pointwise MDP for the Nadaraya-Watson estimator, we need the following additionnal assumptions.
\begin{description}
\item (A3) For any $u \in \mathbb{R}^q$, $t\mapsto \displaystyle \int_{\mathbb{R}^q}\langle u,y\rangle^2f(t,y)dy$ and $t\mapsto \displaystyle \int_{\mathbb{R}^q}\langle u,y\rangle f(t,y)dy$ are continuous at $x$ and $g(x)\neq 0$.
\item (A4) $\lim_{n \to \infty}v_n=\infty$ and $\lim_{n \to \infty}\dfrac{nh_n^d}{v_n^2} =\infty$.
\item (A5) i) There exists an integer $p \geq 2$ such that $\forall s \in \{1,\ldots,p-1\},\forall j \in \{1,\ldots,d\}$,\\ 
$\displaystyle\int_{\mathbb{R}^d}y_j^sK(y)dy_j=0$, and $\displaystyle\int_{\mathbb{R}^d}\big|y_j^pK(y)\big|dy<\infty$.\\
          ii) $\lim_{n \to \infty}v_nh_n^p=0$.\\ 
         iii) $m$ and $g$ are $p$-times differentiable on $\mathbb{R}^d$, and their differentials of order $p$ are bounded and continuous at $x$.
\end{description}
We can now state the MDP for the Nadaraya-Watson estimator.
\begin{theo}\label{mdpregression}\textbf{(Pointwise MDP for the Nadaraya-Watson kernel estimator of the regression)}$ $ \\
Assume that (A1)-(A5) hold. Then, the sequence $\left(v_n\left(r_n(x)-r(x)\right)\right)$ satisfies a LDP with speed $\left(\dfrac{nh_n^d}{v_n^2}\right)$ and good rate function $G_x$ defined in \eqref{rateregression}.
\end{theo}

\subsection{Large and moderate deviations principles for the semi-recursive estimator}
For $a\in ]0,1/d[$, let $\tilde\Psi_{a,x}:\mathbb{R}^q\times\mathbb{R}\rightarrow\mathbb{R}$ and $\tilde{I}_{a,x}:\mathbb{R}^q\times\mathbb{R}\rightarrow\mathbb{R}$ be the functions defined as follows: 
\begin{eqnarray}
\label{psiprimregression}
\tilde\Psi_{a,x}(u,v) & = & \int_{[0,1]\times\mathbb{R}^d\times \mathbb{R}^q}s^{-ad}\left(e^{s^{ad}(\langle u,y\rangle+v)K(z)}-1\right)f(x,y)dsdzdy,\\
\tilde I_{a,x}(t_1,t_2)  &= & \sup_{(u,v) \in \mathbb{R}^q\times\mathbb{R}}\left\{\langle u,t_1\rangle+vt_2-\tilde \Psi_{a,x}(u,v)\right\}.
\label{tauxprimregression}
\end{eqnarray}
 Moreover, let $\tilde J_a$ and $\tilde J^*_a$ be defined as follows: for any $s \in \mathbb{R}^q$,
 \begin{eqnarray}
\label{tauxsemirec1}
\tilde J^*_a(s) & = & \inf_{t \in \mathbb{R}^*}\tilde I_{a,x}(st,t) \\
\label{tauxsemirec2}
\tilde J_a(s) & = & \tilde J^*_a(s)\wedge \tilde I_{a,x}(\vec{0},0).
\end{eqnarray}
Let us give the following additionnal hypotheses.
\begin{description}
\item (A'1) For any $u \in \mathbb{R}^q$, $t\mapsto\displaystyle \int_{\mathbb{R}^q}e^{\alpha\langle u,y\rangle }f(t,y)dy$ is continuous at $x$ uniformly with respect to $\alpha \in[0,1]$.
\end{description}
Condition (C) above is substituted by the following one, 
\begin{description}
\item (C') $\inf_{s \in \mathbb{R}^q}\tilde I_{a,x}(s,0)=\tilde I_{a,x}(\vec{0},0)$. 
\end{description}
Examples for which Condition (C') holds are Examples $1$ and $2$ given for (C).
The following proposition gives the properties of the function $\tilde J _a$.
\begin{prop}\label{barcaprim}
Assume that (A1), (A2), (A'1) and (C') hold. Then,
\begin{description}
\item (i) $\tilde J_a$ is a rate function on $\mathbb{R}^q$. More precisely, for $\alpha\in \mathbb{R}$,
\begin{itemize}
\item if $\alpha<\tilde I_{a,x}(\vec{0},0)$, then $\left\{\tilde J_a(s)\leq \alpha\right \}$ is compact.
\item if $\alpha\geq \tilde I_{a,x}(\vec{0},0)$, then $\left\{\tilde J_a(s)\leq \alpha\right \}=\mathbb{R}^q$.
\end{itemize}
\item (ii) If $\tilde I_{a,x}(\vec{0},0)=\infty$, then $\tilde J_a$ is a good rate function on $\mathbb{R}^q$ and $\tilde J_a=\tilde J^*_a$.
\item (iii) If $\tilde J^*_a(s)<\infty$, then $\tilde J_a(s)=\tilde J^*_a(s)$.
\item (iv) If $\alpha<\tilde I_{a,x}(\vec{0},0)$, then $\left\{\tilde J^*_a(s)\leq \alpha\right \}=\left\{ \tilde J_a(s)\leq \alpha\right \}$. 

\end{description}

\end{prop}
Notice that, like for $J$ and $J^*$, we have
\begin{eqnarray*}
\tilde J_a(s) & = & 
\left\{
\begin{array}{ll}
\tilde J^*_a(s)\ \ \mbox{if} \ \ \tilde J^*_a(s)<\infty\\
\tilde I_{a,x}(\vec{0},0)\ \ \mbox{if} \ \ \tilde J^*_a(s)=\infty.
\end{array}
\right.
\end{eqnarray*}
We can now state the LDP for the semi-recursive kernel estimator of the regression.
\begin{theo}\label{ldpregresprim}\textbf{(Pointwise LDP for the semi-recursive estimator of the regression)}$ $ \\
Set $(h_n)=(cn^{-a})$ with $c>0$ and $0<a<1/d$, and let (A1), (A2), (A'1) and (C') hold. Then, for any open subset $U$ of $\mathbb{R}^q$,  
\begin{eqnarray*}
 \liminf_{n \to \infty}\frac{1}{nh_n^d} \log\mathbb{P}\left[\tilde r_n(x)\in U\right] \geq -\inf_{s\in U}\tilde J^*_a(s), 
\end{eqnarray*}
and for any closed subset $V$ of $\mathbb{R}^q$,
\begin{eqnarray*}
 \limsup_{n \to \infty}\frac{1}{nh_n^d}\log \mathbb{P}\left[\tilde r_n(x)\in V\right] \leq -\inf_{s\in V}\tilde J_a(s).
\end{eqnarray*}
\end{theo}
The comments made for Theorem \ref{ldpregression} are valid for Theorem \ref{ldpregresprim}. In particular, we have the following corollary.
\begin{cor}\label{nonnegativeprim}
Let the assumptions of Theorem \ref{ldpregresprim} hold. If $K$ is a nonnegative kernel such that $\lambda\left(\left\{x\in\mathbb{R}^d, K(x)>0\right\}\right)=\infty$, then the sequence $(\tilde r_n)$ satisfies a LDP with speed $(nh_n^d)$ and good rate function $\tilde J_a$.
\end{cor}

Before stating pointwise MDP for the semi-recursive estimator of the regression, let us recall that a sequence $(u_n)$ is said to vary regularly with exponent $\alpha$ if there exists a function $u$ which varies regularly with exponent $\alpha$ and such that $u_n=u(n)$ for all $n$ (see, for example, Feller \cite{fellerr} page 275).
We will use in the sequel the following property (see Bingham et al. \cite{bingham} page 26). If $(h_n)$ varies regularly with exponent $(-a)$ and if  $\beta a <1$, then  
\begin{eqnarray}
\label{relregress}
\lim_{n \to \infty}\frac{1}{nh_n^\beta}\sum_{i=1}^nh_i^{\beta}= \frac{1}{1-a\beta}.
\end{eqnarray}
We also consider the following condition. 
\begin{eqnarray}
\label{sup}
\sup_{i\leq n}\frac{h_n}{h_i} <\infty.
\end{eqnarray}
(For example, this condition holds when $h_n$ is nonincreasing).
\begin{theo}\label{mdprimregression}\textbf{(Pointwise MDP for the semi-recursive kernel estimator of the regression)}$ $ \\
Assume that $(h_n)$ varies regularly with exponent $(-a)$ with $a\in]0,1/d[$, and satisfies \eqref{sup}. Let (A1)-(A5) hold. Then, the sequence $\left(v_n\left(\tilde r_n(x)-r(x)\right)\right)$ satisfies a LDP with speed $\left(\dfrac{nh_n^d}{v_n^2}\right)$ and good rate function $\tilde G_{a,x}$ defined in \eqref{rateprimregression}.
\end{theo}

\section{Proofs}
The proofs of the results for the Nadaraya-Watson kernel estimator are in many cases similar to those of the semi-recursive kernel estimator of the regression, so we omit some details of the proofs for this last one. \\

First, let us state the following propositions which give the properties of the functions $\Psi_x$, $\tilde \Psi_{a,x}$, $I_x$ and $\tilde I_{a,x}$. Set
\begin{eqnarray*}
S_+=\left\{x\in\mathbb{R}^d,\  K(x)>0\right\} \ \ \mbox{and} \ \ S_-=\left\{x\in\mathbb{R}^d, \ K(x)<0\right\}.
\end{eqnarray*}

\begin{prop}\label{convex1} \textbf{(Properties of \boldmath$\Psi_x$ and $I_x$)}$ $ \\
Let Assumptions (A1) and (A2) hold. Then,
\begin{description} 
\item i) $\Psi_x$ is strictly convex, continuously differentiable on $\mathbb{R}^q\times\mathbb{R}$, and $I_x$ is a good rate function on $\mathbb{R}^q\times\mathbb{R}$.
\item ii) $\nabla \Psi_x$ is an open map and the range of $\Psi_x$ is $\overset{\circ}{\mathcal{D}} (I_x)$. $I_x$ is strictly convex on $\overset{\circ}{\mathcal{D}} (I_x)$ and  for any $t\in \overset{\circ}{\mathcal{D}} (I_x) \subset \mathbb{R}^q\times \mathbb{R}$,
\begin{eqnarray}
\label{maxregression}
I_x(t)=\langle(\nabla \Psi_x)^{-1}(t),t\rangle-\Psi_x\left((\nabla\Psi_x)^{-1}(t)\right).
\end{eqnarray}
\item iii) If $\lambda(S_-)=0$, then $I_x(\vec{0},0)=g(x)\lambda(S_+)$, and for any $t_1\neq \vec{0}$, $I_x(t_1,0)=+\infty$. 
\end{description} 
\end{prop}
\begin{prop}\label{convexprim} \textbf{(Properties of \boldmath  $\tilde\Psi_{a,x}$ and $\tilde I_{a,x}$)}$ $ \\
Let Assumptions (A1) and (A2) hold. Then,
\begin{description} 
\item i) $\tilde\Psi_{a,x}$ is strictly convex, continuously differentiable on $\mathbb{R}^q\times\mathbb{R}$, and $\tilde I_{a,x}$ is a good rate function on $\mathbb{R}^q\times\mathbb{R}$.
\item ii) $\nabla \tilde \Psi_{a,x}$ is an open map and the range of $\tilde\Psi_{a,x}$ is $\overset{\circ}{\mathcal{D}} (\tilde I_{a,x})$. $\tilde I_{a,x}$ is strictly convex on $\overset{\circ}{\mathcal{D}} (\tilde I_{a,x})$, and  for any $t\in \overset{\circ}{\mathcal{D}} (\tilde I_{a,x}) \subset \mathbb{R}^q\times \mathbb{R}$,
\begin{eqnarray}
\label{maxprimregression}
\tilde I_{a,x}(t)=\langle(\nabla \tilde\Psi_{a,x})^{-1}(t),t\rangle-\tilde\Psi_{a,x}\left((\nabla \tilde\Psi_{a,x})^{-1}(t)\right).
\end{eqnarray}
\item iii) If $\lambda(S_-)=0$, then $\tilde I_{a,x}(\vec{0},0)=g(x)\lambda(S_+)/(1-ad)$, and for any $t_1\neq \vec{0}$, $\tilde I_{a,x}(t_1,0)=+\infty$. 
\end{description} 
\end{prop}
 The two following lemmas are used for the proofs of Theorems \ref{ldpregression} and \ref{ldpregresprim}.
\begin{lem}\label{ldpregressioncouple}\textbf{(Pointwise LDP for the sequence \boldmath  $\left(m_n(x),g_n(x)\right)$)}$ $ \\
Let Assumptions (A1) and (A2) hold. Then, the sequence $\left(m_n(x),g_n(x)\right)$ satisfies a LDP with speed $\left(nh_n^d\right)$ and rate function $I_x$ defined in \eqref{tauxregression}.
\end{lem}
\begin{lem}\label{ldprimregressioncouple}\textbf{(Pointwise LDP for the sequence \boldmath  $\left(\tilde m_n(x),\tilde g_n(x)\right)$)}$ $ \\
Set $h_n=cn^{-a}$ with $c>0$ and $a \in ]0,1/d[$, and let Assumptions (A1), (A2) and (A'1) hold. Then, the sequence $\left(\tilde m_n(x),\tilde g_n(x)\right)$ satisfies a LDP with speed $\left(nh_n^d\right)$ and rate function $\tilde I_{a,x}$ defined in \eqref{tauxprimregression}.
\end{lem}

Our proofs are now organized as follows. Lemmas \ref{ldpregressioncouple} and \ref{ldprimregressioncouple} are proved in Section 3.1, Theorems \ref{ldpregression} and \ref{ldpregresprim} in Section 3.2, Theorem \ref{mdpregression} in Section 3.3, Theorem \ref{mdprimregression}  is proved in Section 3.4. Section 3.5 is devoted to the proof of Propositions \ref{convex1} and \ref{convexprim} on the rate functions $I_x$ and $I_{a,x}$. Propositions \ref{barca} and \ref{barcaprim} are proved in Section 3.6.

\subsection{Proof of Lemmas \ref{ldpregressioncouple} and \ref{ldprimregressioncouple}} 
\subsubsection{Proof of Lemma \ref{ldpregressioncouple}}
For any $w=(u,v) \in \mathbb{R}^q\times \mathbb{R}$, set
\begin{eqnarray*}
\Psi_n(x) & =& \left(m_n(x),g_n(x)\right),\\
\Lambda_{n,x}(w) & = & \frac{1}{nh_n^d}\log\mathbb{E}\left[\exp\left(nh_n^d\langle w,\Psi_n(x)\rangle\right)\right].
\end{eqnarray*}
Let us at first assume that the following lemma holds.
\begin{lem}\label{cum} \textbf{(Convergence of \boldmath $\Lambda_{n,x}$)} $ $\\
Assume that (A1) and (A2) hold, then
\begin{eqnarray}
\label{convregression}
\lim_{n \to \infty}\Lambda_{n,x}(u,v) & =& \Psi_x(u,v),
\end{eqnarray}
where $\Psi_x$ is defined in \eqref{psiregression}.
\end{lem}
To prove Lemma \ref{ldpregressioncouple}, we apply Proposition \ref{convex1}, Lemma \ref{cum} and the G\"artner-Ellis Theorem (see Dembo and Zeitouni \cite{dembzeitr}). 
Proposition \ref{convex1} ensures that $\Psi_x$ is essentially smooth, lower semicontinuous function so that Lemma \ref{ldpregressioncouple} follows from the G\"artner-Ellis Theorem.\\
Let us now prove Lemma \ref{cum}.
Set
\begin{eqnarray*}
Z_i=\left[\langle u,Y_i\rangle+v\right]K\left(\frac{x-X_i}{h_n}\right).
\end{eqnarray*}
For any $\left(u,v\right) \in \mathbb{R}^q\times \mathbb{R}$, we have 
\begin{eqnarray*}
\Lambda_{n,x}(u,v) & = & \frac{1}{nh_n^d}\log\mathbb{E}\left[\exp\left(\sum_{i=1}^nZ_i\right)\right],
\end{eqnarray*}
and, since the random vectors $(X_i,Y_i)$, $i=1,\dots,n$ are independent and identically distributed, we get
\begin{eqnarray*}
\Lambda_{n,x}(u,v) & = & \frac{1}{h_n^d}\log\mathbb{E}\left[e^{Z_n}\right].
\end{eqnarray*}
A Taylor's expansion implies that there exists $c_n$ between $1$ and $\mathbb{E}\left[e^{Z_n}\right]$ such that
\begin{eqnarray*}
\Lambda_{n,x}(u,v)
& = & 
\frac{1}{h_n^d}\mathbb{E}\left[e^{Z_n}-1\right]-\frac{1}{2c_n^2h_n^d}\left(\mathbb{E}\left[e^{Z_n}-1\right]\right)^2\\
& = &  
\frac{1}{h_n^d}\int_{\mathbb{R}^d \times \mathbb{R}^q}
\left[e^{\left(\langle u,y\rangle+v\right)K\left(\frac{x-s}{h_n}\right)}-1\right]f(s,y)dsdy-R_{n,x}^{(1)}(u,v) \\ & = & 
\Psi_x(u,v)-R_{n,x}^{(1)}(u,v)+R_{n,x}^{(2)}(u,v),
\end{eqnarray*}
with 
\begin{eqnarray*}
R_{n,x}^{(1)}(u,v) & = & \frac{1}{2c_n^2h_n^d}\left(\mathbb{E}\left[e^{Z_n}-1\right]\right)^2,\\
R_{n,x}^{(2)}(u,v) & = &  \int_{\mathbb{R}^d\times\mathbb{R}^q}\left[e^{\left(\langle u,y\rangle+v\right)K\left(z\right)}-1\right]\left[f(x-h_nz,y)-f(x,y)\right]dzdy.
\end{eqnarray*}
Let us prove that 
\begin{eqnarray}
\label{rn2}
\lim_{n \to \infty}R_{n,x}^{(2)}(u,v) &= & 0.
\end{eqnarray}
Set $A>0$ and $\epsilon>0$; we then have
\begin{eqnarray}
\label{qua}
R_{n,x}^{(2)}(u,v) & = & \int_{\left\{\|z\|\leq A\right\}\times\mathbb{R}^q}\left[e^{\left(\langle u,y\rangle+v\right)K\left(z\right)}-1\right]\left[f(x-h_nz,y)-f(x,y)\right]dzdy   \nonumber \\ \mbox{} & &+\int_{\left\{\|z\| > A\right\}\times\mathbb{R}^q}\left[e^{\left(\langle u,y\rangle+v\right)K\left(z\right)}-1\right]\left[f(x-h_nz,y)-f(x,y)\right]dzdy.
\end{eqnarray}
Next, since for any $t \in \mathbb{R}$, $\left|e^t-1\right| \leq |t|e^{|t|}$, we have
\begin{eqnarray}
\label{quo}
\lefteqn{\int_{\left\{\|z\| > A\right\}\times\mathbb{R}^q}\left|e^{\left(\langle u,y\rangle+v\right)K\left(z\right)}-1\right|\left|f(x-h_nz,y)-f(x,y)\right|dzdy} \nonumber\\ &\leq & 
\int_{\left\{\|z\| > A\right\}\times\mathbb{R}^q}\left|\langle u,y\rangle+v\right|\left|K\left(z\right)\right|e^{\left|\langle u,y\rangle+v\right|\left|K\left(z\right)\right|}\left|f(x-h_nz,y)-f(x,y)\right|dzdy \nonumber\\ & \leq & 
\int_{\left\{\|z\| > A\right\}\times\mathbb{R}^q}\left|\langle u,y\rangle+v\right|\left|K\left(z\right)\right|e^{\left|\langle u,y\rangle+v\right|\left|K\left(z\right)\right|}f(x-h_nz,y)dzdy \nonumber\\ \mbox{} & & +\int_{\left\{\|z\| > A\right\}\times\mathbb{R}^q}\left|\langle u,y\rangle+v\right|\left|K\left(z\right)\right|e^{\left|\langle u,y\rangle+v\right|\left|K\left(z\right)\right|}f(x,y)dzdy  \nonumber\\ & \leq &
e^{|v|\|K\|_{\infty}}\|u\|\int_{\left\{\|z\| > A\right\}}\left|K\left(z\right)\right|\left[\int_{\mathbb{R}^q}\|y\|e^{\|K\|_{\infty}\|u\|\|y\|}f(x-h_nz,y)dy\right]dz \nonumber \\ \mbox{} & & + e^{|v|\|K\|_{\infty}}\left|v\right| \int_{\left\{\|z\| > A\right\}}\left|K\left(z\right)\right|\left[\int_{\mathbb{R}^q}e^{\|K\|_{\infty}\|u\|\|y\|}f(x-h_nz,y)dy\right]dz\nonumber\\ \mbox{} & & +
e^{|v|\|K\|_{\infty}}\|u\|\int_{\left\{\|z\| > A\right\}}\left|K\left(z\right)\right|dz\int_{\mathbb{R}^q}\|y\|e^{\|K\|_{\infty}\|u\|\|y\|}f(x,y)dy \nonumber \\ \mbox{} & & + e^{|v|\|K\|_{\infty}}\left|v\right| \int_{\left\{\|z\| > A\right\}}\left|K\left(z\right)\right|dz\int_{\mathbb{R}^q}e^{\|K\|_{\infty}\| u\|\|y\|}f(x,y)dy \nonumber\\ & \leq & 
B\int_{\left\{\|z\| > A\right\}}\left|K\left(z\right)\right|dz,
\end{eqnarray}
where $B$ is a constant ; this last inequality follows from \eqref{borne} and from the fact that $K$ is bounded. Now, since $K$ is integrable, we can choose $A$ such that
\begin{eqnarray}
\label{secondterm}
\int_{\left\{\|z\| > A\right\}\times\mathbb{R}^q}\left|e^{\left(\langle u,y\rangle+v\right)K\left(z\right)}-1\right|\left|f(x-h_nz,y)-f(x,y)\right|dzdy \leq \frac{\epsilon}{2}.
\end{eqnarray}
Now, observe that
\begin{eqnarray}
\label{xon1}
\lefteqn{\int_{\left\{\|z\|\leq A\right\}\times\mathbb{R}^q}\left[e^{\left(\langle u,y\rangle+v\right)K\left(z\right)}-1\right]\left[f(x-h_nz,y)-f(x,y)\right]dzdy} \nonumber\\ & = & \int_{\left\{\|z\|\leq A\right\}}e^{vK(z)}\left[\int_{\mathbb{R}^q}e^{\langle u,y \rangle K(z)}\left(f(x-h_nz,y)-f(x,y)\right)dy\right]dz \\ \mbox{} & & -
\int_{\left\{\|z\|\leq A\right\}}\left[\int_{\mathbb{R}^q}\left(f(x-h_nz,y)-f(x,y)\right)dy\right]dz.
\label{xon2}
\end{eqnarray}
Assumption (A2) together with \eqref{borne}, and the dominated convergence theorem ensure that both integrals in \eqref{xon1} and \eqref{xon2} converge to $0$. We deduce that for $n$ large enough,
\begin{eqnarray}
\label{firsterm}
\left|\int_{\left\{\|z\|\leq A\right\}\times\mathbb{R}^q}\left[e^{\left(\langle u,y\rangle+v\right)K\left(z\right)}-1\right]\left[f(x-h_nz,y)-f(x,y)\right]dzdy\right| \leq \frac{\epsilon}{2},
\end{eqnarray}
so that \eqref{rn2} follows from \eqref{secondterm} and \eqref{firsterm}.\\ 
Let us now consider $R_{n,x}^{(1)}$ ; since $c_n$ is between $1$ and $\mathbb{E}\left[e^{Z_n}\right]$, we get
\begin{eqnarray*}
\frac{1}{c_n} & \leq & \max\left\{1,\frac{1}{\mathbb{E}\left(e^{Z_n}\right)}\right\}. 
\end{eqnarray*}
By Jensen's inequality, we obtain
\begin{eqnarray*}
\frac{1}{\mathbb{E}\left[e^{Z_n}\right]} \leq \frac{1}{e^{\mathbb{E}\left[Z_n\right]}}.
\end{eqnarray*}
Observe that
\begin{eqnarray*}
\left|\mathbb{E}\left(Z_n\right)\right| & = &\left| \mathbb{E}\left[\left(\langle u,Y_n\rangle+v\right)K\left(\frac{x-X_n}{h_n}\right)\right]\right| \\ 
& \leq & 
\int_{\mathbb{R}^d \times \mathbb{R}^q}\|u\|\|y\|\left| K\left(\frac{x-s}{h_n}\right)\right|f(s,y)dsdy+|v| \int_{\mathbb{R}^d \times \mathbb{R}^q}\left|K\left(\frac{x-s}{h_n}\right)\right|f(s,y)dsdy \\ 
& \leq & 
h_n^d\left(\|u\|\int_{\mathbb{R}^d }\left|K(z)\right|\left[\int_{\mathbb{R}^q}\|y\|f(x-h_nz,y)dy\right]dz+|v| \int_{\mathbb{R}^d }\left|K(z)\right|\left[\int_{\mathbb{R}^q}f(x-h_nz,y)dy \right]dz\right),
\end{eqnarray*}
which goes to $0$ in view of \eqref{borne} and since $\lim_{n \to \infty}h_n=0$. We deduce that there exists $c \in \mathbb{R}^*_+$ such that 
\begin{eqnarray*}
\frac{1}{c_n^2} \leq c.
\end{eqnarray*}
Noting that by \eqref{borne},
\begin{eqnarray*}
\mathbb{E}\left|e^{Z_n}-1\right| & \leq & h_n^d\int_{\mathbb{R}^d\times\mathbb{R}^q}\left|\langle u,y\rangle+v\right||K(z)|e^{\left|\langle u,y\rangle+v\right||K(z)|}f(x-h_nz,y)dydz\\ & \leq & 
Bh_n^d,
\end{eqnarray*}
where $B$ is a constant.
It follows that 
\begin{eqnarray*}
\lim_{n \to \infty}R_{n,x}^{(1)}(u,v)=0,
\end{eqnarray*}
which proves Lemma \ref{cum}. $\blacksquare$

\subsubsection{Proof of Lemma \ref{ldprimregressioncouple}}
Similarly as the proof of Lemma \ref{ldpregressioncouple},
for any $w=(u,v) \in \mathbb{R}^q\times \mathbb{R}$, set
\begin{eqnarray*}
\tilde \Psi_n(x) & =& \left(\tilde m_n(x),\tilde g_n(x)\right),\\
\tilde \Lambda_{n,x}(w) & = & \frac{1}{nh_n^d}\log\mathbb{E}\left[\exp\left(nh_n^d\langle w,\tilde \Psi_n(x)\rangle\right)\right].
\end{eqnarray*}
When $h_n=cn^{-a}$, $c>0$ and $0<a<1/d$, assume for the moment that 
\begin{eqnarray}
\label{secondcumulant}
\lim_{n \to \infty}\tilde \Lambda_{n,x}(u,w) & = & \tilde \Psi_{a,x}(u,v),
\end{eqnarray}
where $\tilde \Psi_{a,x}$ is defined in \eqref{psiprimregression}. The conclusion of Lemma \ref{ldprimregressioncouple} follows from Proposition \ref{convexprim} and again the G\"artner-Ellis Theorem. \\

Let us now prove \eqref{secondcumulant}. Set
\begin{eqnarray*}
M_i=\left[\langle u,Y_i\rangle+v\right]K\left(\frac{x-X_i}{h_i}\right),
\end{eqnarray*}
then, for $\left(u,v\right) \in \mathbb{R}^q\times \mathbb{R}$, 
\begin{eqnarray*}
\tilde\Lambda_{n,x}(u,v) & = & \frac{1}{nh_n^d}\log\mathbb{E}\left[\exp\left(\sum_{i=1}^nM_i\frac{h_n^d}{h_i^d}\right)\right]\\ & = & 
\frac{1}{nh_n^d}\sum_{i=1}^n\log\mathbb{E}\left[\exp\left(M_i\frac{h_n^d}{h_i^d}\right)\right].
\end{eqnarray*}
By Taylor expansion, there exists $b_{i,n}$ between $1$ and $\mathbb{E}\left[\exp\left(M_i\frac{h_n^d}{h_i^d}\right)\right]$ such that 
\begin{eqnarray*}
\log\mathbb{E}\left[\exp\left(M_i\frac{h_n^d}{h_i^d}\right)\right] &= & \mathbb{E}\left[\exp\left(M_i\frac{h_n^d}{h_i^d}\right)-1\right]-\frac{1}{2b_{i,n}^2}\left(\mathbb{E}\left[\exp\left(M_i\frac{h_n^d}{h_i^d}\right)-1\right]\right)^2.
\end{eqnarray*}
Noting that $h_n=cn^{-a}$ with $c>0$ and $a\in]0,1/d[$, $\tilde\Lambda_{n,x}$ can be rewritten as
\begin{eqnarray*}
\tilde\Lambda_{n,x}(u,v) & = & \frac{1}{h_n^d}\sum_{i=1}^n \mathbb{E}\left[\exp\left(M_i\frac{h_n^d}{h_i^d}\right)-1\right]-\frac{1}{2nh_n^d}\sum_{i=1}^n\frac{1}{b_{i,n}^2}\left(\mathbb{E}\left[\exp\left(M_i\frac{h_n^d}{h_i^d}\right)-1\right]\right)^2\\ &= &
\frac{1}{n}\sum_{i=1}^n\left(\frac{i}{n}\right)^{-ad}\int_{\mathbb{R}^d\times\mathbb{R}^q}\left[e^{\left(\frac{i}{n}\right)^{ad}\left(\langle u,y\rangle+v\right)K\left(z\right)}-1\right]f(x,y)dzdy -R_{n,x}^{(1)}(u,v)+R_{n,x}^{(2)}(u,v),
\end{eqnarray*}
with
\begin{eqnarray*}
\tilde R_{n,x}^{(1)}(u,v) &= &\frac{1}{2nh_n^d}\sum_{i=1}^n\frac{1}{b_{i,n}^2}\left(\mathbb{E}\left[\exp\left(M_i\frac{h_n^d}{h_i^d}\right)-1\right]\right)^2\\ 
\tilde R_{n,x}^{(2)}(u,v) & = & 
 \frac{1}{nh_n^d}\sum_{i=1}^nh_i^d\int_{\mathbb{R}^d\times\mathbb{R}^q}\left[e^{\frac{h_n^d}{h_i^d}\left(\langle u,y\rangle+v\right)K\left(z\right)}-1\right]\left[f(x-h_iz,y)-f(x,y)\right]dzdy.
\end{eqnarray*}

Since $b_{i,n}$ is between $1$ and $\mathbb{E}\left[\exp\left(M_i\frac{h_n^d}{h_i^d}\right)\right]$, we have
\begin{eqnarray*}
\frac{1}{b_{i,n}} & \leq & \max\left\{1,\frac{1}{\mathbb{E}\left(e^{M_i\frac{h_n^d}{h_i^d}}\right)}\right\}. 
\end{eqnarray*}
By Jensen's inequality, we obtain
\begin{eqnarray*}
\frac{1}{\mathbb{E}\left[e^{M_i\frac{h_n^d}{h_i^d}}\right]} \leq \frac{1}{e^{\mathbb{E}\left[M_i\frac{h_n^d}{h_i^d}\right]}}.
\end{eqnarray*}
Observe that 
\begin{eqnarray*}
\lefteqn{\left|\mathbb{E}\left(M_i\frac{h_n^d}{h_i^d}\right)\right|}\\ &\leq& \frac{h_n^d}{h_i^d}\int_{\mathbb{R}^d}\left|\langle u,y\rangle+v\right|\left|K\left(\frac{x-z}{h_i}\right)\right|f(z,y)dzdy \\ & \leq &  
h_n^d\left(\|u\|\int_{\mathbb{R}^d\times\mathbb{R}^q }\left|K(z)\right|\left[\int_{\mathbb{R}^q}\|y\|f(x-h_iz,y)dy\right]dz +|v| \int_{\mathbb{R}^d }\left|K(z)\right|\left[\int_{\mathbb{R}^q}f(x-h_iz,y)dy \right]dz\right),
\end{eqnarray*}
which goes to $0$ in view of \eqref{borne} and since $\lim_{n \to \infty}h_n=0$. We deduce that the sequence $\left(\mathbb{E}\left[\exp\left(M_i\frac{h_n^d}{h_i^d}\right)\right]\right)$ is bounded, so that
there exists $c>0$ such that 
\begin{eqnarray*}
\frac{1}{b_{i,n}^2} \leq c,
\end{eqnarray*}
and thus
\begin{eqnarray*}
\tilde R_{n,x}^{(1)}(u,v) & \leq & \frac{c}{2nh_n^d}\sum_{i=1}^n\left(\mathbb{E}\left[e^{M_i\frac{h_n^d}{h_i^d}}-1\right]\right)^2.
\end{eqnarray*}
Now, in view of \eqref{borne}, and since $K$ is bounded integrable, we have
\begin{eqnarray*} 
\mathbb{E}\left|e^{M_i\frac{h_n^d}{h_i^d}}-1 \right| & \leq & \mathbb{E}\left[ \left|M_i\frac{h_n^d}{h_i^d}\right|e^{\left|M_i\frac{h_n^d}{h_i^d}\right|}\right]\\ & \leq & 
\frac{h_n^d}{h_i^d}\int_{\mathbb{R}^d}\left|\langle u,y\rangle+v\right|\left|K\left( \frac{x-s}{h_i}\right)\right|e^{\left|\langle u,y\rangle+v\right|\left|K\left(\frac{x-s}{h_i}\right)\right|}f(s,y)dsdy\\ & \leq & h_n^d\int_{\mathbb{R}^d \times \mathbb{R}^q}\left|\langle u,y\rangle+v\right|\left|K\left(z\right)\right|e^{\left|\langle u,y\rangle+v\right|\left|K\left(z\right)\right|}f(x-h_iz,y)dzdy \\ & \leq & Bh_n^d,
\end{eqnarray*}
where $B$ is a constant. Thus
\begin{eqnarray*}
\left|\tilde R_{n,x}^{(1)}(u,v)\right| & \leq & \frac{cB^2}{2}h_n^{d},
\end{eqnarray*}
and 
\begin{eqnarray*}
\lim_{n \to \infty}\left|\tilde R_{n,x}^{(1)}(u,v)\right|=0.
\end{eqnarray*}
Let us now consider $R_{n,x}^{(2)}$.
Set $A>0$ and $\epsilon>0$ ; we then have
\begin{eqnarray*}
\lefteqn{\tilde R_{n,x}^{(2)}(u,v)}\\ & = & 
 \frac{1}{nh_n^d}\sum_{i=1}^nh_i^d\int_{\left\{\|z\|\leq A\right\}\times\mathbb{R}^q}\left[e^{\frac{h_n^d}{h_i^d}\left(\langle u,y\rangle+v\right)K\left(z\right)}-1\right]\left[f(x-h_iz,y)-f(x,y)\right]dzdy \\ \mbox{} && +
 \frac{1}{nh_n^d}\sum_{i=1}^nh_i^d\int_{\left\{\|z\|> A\right\}\times\mathbb{R}^q}\left[e^{\frac{h_n^d}{h_i^d}\left(\langle u,y\rangle+v\right)K\left(z\right)}-1\right]\left[f(x-h_iz,y)-f(x,y)\right]dzdy \\& = &  I+II .
\end{eqnarray*}
Since $|e^t-1| \leq |t|e^{|t|}$, it follows that
\begin{eqnarray*}
\left| II\right| & \leq & \frac{1}{n}\sum_{i=1}^n\int_{\left\{\|z\|> A\right\}\times\mathbb{R}^q}\left|\langle u,y\rangle+v\right||K\left(z\right)|e^{\left|\langle u,y\rangle+v\right||K\left(z\right)|}\left|f(x-h_iz,y)-f(x,y)\right|dzdy.
\end{eqnarray*}
Using the same argument as in \eqref{quo}, it holds that
\begin{eqnarray*}
\left|II\right| \leq \frac{\epsilon}{2}.
\end{eqnarray*}
Now, for $I$, we write
\begin{eqnarray*}
I& = & \frac{1}{nh_n^d}\sum_{i=1}^nh_i^d\int_{\left\{\|z\|\leq A\right\}\times\mathbb{R}^q}e^{\frac{h_n^d}{h_i^d}\left(\langle u,y\rangle+v\right)K\left(z\right)}\left[f(x-h_iz,y)-f(x,y)\right]dzdy \\ \mbox{} & &- \frac{1}{nh_n^d}\sum_{i=1}^nh_i^d\int_{\left\{\|z\|\leq A\right\}\times\mathbb{R}^q}\left[f(x-h_iz,y)-f(x,y)\right]dzdy .
\end{eqnarray*}
On the one hand, Assumption (A2) with $u=0$ ensures that
\begin{eqnarray*}
\lim_{i \to \infty}\int_{\left\{\|z\|\leq A\right\}\times\mathbb{R}^q}\left[f(x-h_iz,y)-f(x,y)\right]dzdy &= & 0.
\end{eqnarray*}
Moreover, since $ad<1$, \eqref{relregress} ensures that
\begin{eqnarray*}
\lim_{n \to \infty}\frac{1}{nh_n^d}\sum_{i=1}^nh_i^d\int_{\left\{\|z\|\leq A\right\}\times\mathbb{R}^q}\left[f(x-h_iz,y)-f(x,y)\right]dzdy & = & 0,	
\end{eqnarray*}
so that for $n$ large enough,
\begin{eqnarray}
\label{eps1}
\left|\frac{1}{nh_n^d}\sum_{i=1}^nh_i^d\int_{\left\{\|z\|\leq A\right\}\times\mathbb{R}^q}\left[f(x-h_iz,y)-f(x,y)\right]dzdy\right|	\leq \frac{\epsilon}{4}.					  
\end{eqnarray}
On the other hand,
since for $i\leq n$, $0\leq \frac{h_n^d}{h_i^d}\leq 1$, by Assumption (A'1), there exists $n_0\in \mathbb{N}$ such that for any $i>n_0$,
\begin{eqnarray*}
\left|\int_{\mathbb{R}^q}e^{\frac{h_n^d}{h_i^d}\left(\langle u,y\rangle+v\right)K\left(z\right)}\left[f(x-h_iz,y)-f(x,y)\right]dy\right| \leq \frac{\epsilon}{8(2-ad)\int_{\|z\| \leq A}dz} \ \ \forall n>i.
\end{eqnarray*}
Noting that by \eqref{borne}, for any $\alpha \in [0,1]$, 
\begin{eqnarray*}
\sup_t\left|\int_{\mathbb{R}^q}e^{\alpha\langle u,y\rangle}f(t,y)dy\right| &\leq& 
\sup_t\int_{\mathbb{R}^q}e^{\|u\|\|y\|}f(t,y)dy <  \infty.
\end{eqnarray*}
Since $ad<1$, by \eqref{relregress}, we get for $n$ sufficiently large
\begin{eqnarray*}
\left|\frac{1}{nh_n^d}\sum_{i=n_0+1}^nh_i^d\int_{\left\{\|z\|\leq A\right\}\times\mathbb{R}^q}e^{\frac{h_n^d}{h_i^d}\left(\langle u,y\rangle+v\right)K\left(z\right)}\left[f(x-h_iz,y)-f(x,y)\right]dydz\right| & \leq & \frac{\epsilon}{8}.
\end{eqnarray*}
Now, for $n$ large enough, in view of \eqref{borne},  
\begin{eqnarray*}
\lefteqn{\left|\frac{1}{nh_n^d}\sum_{i=1}^{n_0}h_i^d\int_{\left\{\|z\|\leq A\right\}\times\mathbb{R}^q}e^{\frac{h_n^d}{h_i^d}\left(\langle u,y\rangle+v\right)K\left(z\right)}\left[f(x-h_iz,y)-f(x,y)\right]dydz\right|}\\ & \leq & 
\frac{1}{nh_n^d}\sum_{i=1}^{n_0}h_i^d\int_{\left\{\|z\|\leq A\right\}\times\mathbb{R}^q}e^{\left|\langle u,y\rangle+v\right||K\left(z\right)|}\left|f(x-h_iz,y)-f(x,y)\right|dydz
\\&\leq &\frac{\epsilon}{8}.
\end{eqnarray*}
It follows that for $n$ large enough, 
\begin{eqnarray}
\label{eps2}
\left|\frac{1}{nh_n^d}\sum_{i=1}^nh_i^d\int_{\left\{\|z\|\leq A\right\}\times\mathbb{R}^q}e^{\frac{h_n^d}{h_i^d}\left(\langle u,y\rangle+v\right)K\left(z\right)}\left[f(x-h_iz,y)-f(x,y)\right]dydz\right| & \leq & \frac{\epsilon}{4}.
\end{eqnarray}
The combination of \eqref{eps1} and \eqref{eps2} ensures that $\left|I\right|\leq \dfrac{\epsilon}{2}$, which ensures that
\begin{eqnarray*}
\lim_{n\to \infty}\left|\tilde R_{n,x}^{(2)}(u)\right| & = & 0.
\end{eqnarray*}
 Hence, \eqref{secondcumulant} follows from analysis considerations. $\blacksquare$

\subsection{Proof of Theorems \ref{ldpregression} and \ref{ldpregresprim}}
Let us consider the following functions defined as:
\begin{eqnarray*}
H_1:\mathbb{R}^q\times \mathbb{R}^* &\rightarrow & \mathbb{R}^q\\
(\alpha,\beta) & \mapsto &  \frac{\alpha}{\beta},
\end{eqnarray*}
and
\begin{eqnarray*}
H_2:\mathbb{R}^q\times \mathbb{R} &\rightarrow & \mathbb{R}^q\\
(\alpha,\beta) & \mapsto &  
\left\{
\begin{array}{ll}
\dfrac{\alpha}{\beta} \ \ \mbox{if} \ \ \beta \neq 0\\ \nonumber
0 \ \ \mbox{otherwise}.
\end{array}
\right.
\end{eqnarray*}
\subsubsection{Proof of Theorem \ref{ldpregression}}
\begin{description}
 \item i) Let $U$ be an open subset of $\mathbb{R}^q$, we have
\begin{eqnarray}
\label{probregression}
\frac{1}{nh_n^d}\log\mathbb{P}\left[r_n(x)\in U\right] & = & \frac{1}{nh_n^d}\log\mathbb{P}\left[\left(m_n(x),g_n(x)\right)\in H_2^{-1}(U)\right].
\end{eqnarray}
Observe that $H_1^{-1}(U)\subset H_2^{-1}(U)$ and $H_1^{-1}(U)$ is an open subset on $\mathbb{R}^q\times\mathbb{R}^*$ which is open, it follows that $H_1^{-1}(U)$ is an open subset on $\mathbb{R}^q\times\mathbb{R}$. We deduce from \eqref{probregression} that
\begin{eqnarray*}
\frac{1}{nh_n^d}\log\mathbb{P}\left[r_n(x)\in U\right] & \geq & \frac{1}{nh_n^d}\log\mathbb{P}\left[\left(m_n(x),g_n(x)\right)\in H_1^{-1}(U)\right].
\end{eqnarray*}
The application of Lemma \ref{ldpregressioncouple} ensures that
\begin{eqnarray*}
\liminf_{n \to \infty}\frac{1}{nh_n^d}\log\mathbb{P}\left[r_n(x)\in U\right] & \geq  & 
\liminf_{n \to \infty}\frac{1}{nh_n^d}\log\mathbb{P}\left[\left(m_n(x),g_n(x)\right)\in H_1^{-1}(U)\right]  \nonumber \\ & \geq & 
-\inf_{\left(t_1,t_2\right) \in H_1^{-1}(U)}I_x(t_1,t_2)= - \inf_{s\in U}J^*(s),
\end{eqnarray*}
and the first part of Theorem \ref{ldpregression} is proved. 

\item ii)  Let $V$ be a closed subset of $\mathbb{R}^q$, we have
\begin{eqnarray*}
\frac{1}{nh_n^d}\log\mathbb{P}\left[r_n(x)\in V\right] & = & \frac{1}{nh_n^d}\log\mathbb{P}\left[\left(m_n(x),g_n(x)\right)\in  H_2^{-1}(V)\right]\\ & \leq & \frac{1}{nh_n^d}\log\mathbb{P}\left[\left(m_n(x),g_n(x)\right)\in \overline{H_2^{-1}(V)}\right].
\end{eqnarray*}
Now, observe that $\overline{H_2^{-1}(V)}=H_1^{-1}(V)\cup A$ where $A\subset \mathbb{R}^q\times\{0\}$ and $\big(\vec{0},0\big)\in A$ (since for  any $s \in \mathbb{R}^q$, $(\vec{0},0)\in\overline{H_2^{-1}(s)}$). The application of Lemma \ref{ldpregressioncouple} again ensures that 
\begin{eqnarray*}
\limsup_{n \to \infty}\frac{1}{nh_n^d}\log\mathbb{P}\left[r_n(x)\in V\right] & \leq & -\inf_{(s,t)\in H_1^{-1}(V)\cup A }I_x(s,t) \\ & \leq &
-\inf_{(s,t)\in H_1^{-1}(V)\cup \{(\vec{0},0)\}}I_x(s,t)\\ & \leq &
-\inf_{s \in V, \ t \in \mathbb{R}}I_x(st,t) \\  & \leq & -\inf_{s \in V, \ t \in \mathbb{R}}\hat{I}_x(s,t)\\ & \leq & 
-\inf_{s \in V}J(s),
\end{eqnarray*}
where the second inequality comes from Condition (C); this concludes the proof of Theorem \ref{ldpregression}. $\blacksquare$\\
\subsubsection{Proof of Theorem \ref{ldpregresprim}}
Applying Lemma \ref{ldprimregressioncouple}, Theorem \ref{ldpregresprim} is proved by following the same approach as for the proof of Theorem \ref{ldpregression} with replacing $m_n$, $g_n$, $J^*$ and $J$ by $\tilde m_n$, $\tilde g_n$, $\tilde J^*_a$ and $\tilde J_a$ respectively. $\blacksquare$

\end{description}

\subsection{Proof of Theorem \ref{mdpregression}}
Set 
\begin{eqnarray*}
B_n(x) & = & \frac{1}{g(x)}\left(m_n(x)-m(x)\right)-\frac{r(x)}{g(x)}\left(g_n(x)-g(x)\right).
\end{eqnarray*}
Let us at first state the two following lemmas.
\begin{lem}\label{mdpbn} $ $
Under the assumptions of Theorem \ref{mdpregression}, the sequence $\left(v_n\left(B_n(x)-\mathbb{E}\left(B_n(x)\right)\right)\right)$ satisfies a LDP with speed $\left(\frac{nh_n^d}{v_n^2}\right)$ and good rate function $G_x$.
\end{lem}
\begin{lem}\label{biaisregression}$ $
Under the assumptions of Theorem \ref{mdpregression},
\begin{eqnarray}
\label{esperanceregression}
\lim_{n \to \infty} v_n\mathbb{E}(B_n(x)) & =&0.
\end{eqnarray}

\end{lem}

We first show that how Theorem \ref{mdpregression} can be deduced from the application of Lemmas \ref{mdpbn} and \ref{biaisregression}, and then prove Lemmas \ref{mdpbn} and \ref{biaisregression} successively.
\subsubsection{Proof of Theorem \ref{mdpregression}}
Lemmas \ref{mdpbn} and \ref{biaisregression} imply that the sequence $\left(v_nB_n(x)\right)$ satisfies a LDP with speed $\left(\dfrac{nh_n^d}{v_n^2}\right)$ and good rate function $G_x$. To prove Theorem \ref{mdpregression}, we show that $\left(v_n(r_n-r)\right)$ and $\left(v_nB_n\right)$ are exponentially contiguous.\\
Let us first note that, for $x$ such that $g_n(x)\neq 0$, we have:
\begin{eqnarray*}
r_n(x)-r(x) & = & \frac{m_n(x)}{g_n(x)}-\frac{m(x)}{g(x)} \\ & = & 
\frac{\left(m_n(x)-m(x)\right)g(x)+\left(g(x)-g_n(x)\right)m(x)}{g_n(x)g(x)}\\ & = & B_n(x)\frac{g(x)}{g_n(x)}.
\end{eqnarray*}
It follows that, for any $\delta>0$, we have
\begin{eqnarray*}
\label{equation1regression}
\lefteqn{\mathbb{P}\Big[v_n\|\left(r_n(x)-r(x)\right)-B_n(x)\|>\delta\Big]}\\ & \leq &
 \mathbb{P}\Big[v_n\|B_n(x)\left(\frac{g(x)}{g_n(x)}-1\right)\|>\delta \ \ \mbox{and} \ \ g_n(x)\neq 0\Big] +\mathbb{P}\left[g_n(x)=0\right] \\ & \leq & 
\mathbb{P}\left[\sqrt{v_n}\|B_n(x)\|>\delta\right]+\mathbb{P}\left[\sqrt{v_n}\left|g(x)-g_n(x)\right| >\delta\left|g_n(x)\right|\right]+\mathbb{P}\left[\left|g(x)-g_n(x)\right| >\frac{g(x)}{2}\right]\\
& \leq & 
\mathbb{P}\left[\sqrt{v_n}\|B_n(x)\|>\delta\right]+\mathbb{P}\left[\sqrt{v_n}\left|g(x)-g_n(x)\right| >\delta\left|g_n(x)\right| \ \ \mbox{and} \ \ \frac{g_n(x)}{g(x)}>\frac{1}{2}\right]\\ \mbox{} && + \mathbb{P}\left[\frac{g_n(x)}{g(x)}\leq\frac{1}{2} \right]+\mathbb{P}\left[\left|g_n(x)-g(x)\right| >\frac{g(x)}{2}\right] \\ 
& \leq & 
\mathbb{P}\left[\sqrt{v_n}\|B_n(x)\|>\delta\right]+\mathbb{P}\left[\sqrt{v_n}\left|g(x)-g_n(x)\right| >\delta\frac{g(x)}{2}\right]+\mathbb{P}\left[g(x)-g_n(x) \geq\frac{g(x)}{2}\right] \\ \mbox{} && +\mathbb{P}\left[\left|g_n(x)-g(x)\right| >\frac{g(x)}{2}\right].
\end{eqnarray*}
 Since $\lim_{n \to \infty}v_n=\infty$, it follows that, for $n$ large enough,
\begin{eqnarray*}
\lefteqn{\mathbb{P}\Big[v_n\|\left(r_n(x)-r(x)\right)-B_n(x)\|>\delta\Big]}\\ & \leq & 
4\max\left\{\mathbb{P}\left[\sqrt{v_n}\|B_n(x)\|>\delta\right] \ ;  \ \mathbb{P}\left[\sqrt{v_n}\left|g(x)-g_n(x)\right| >\delta\frac{g(x)}{2}\right] \right\},
\end{eqnarray*}
and thus
\begin{eqnarray*}
\lefteqn{\frac{v_n^2}{nh_n^d}\log\mathbb{P}\Big[v_n\|\left(r_n(x)-r(x)\right)-B_n(x)\|>\delta\Big]} \\ 
& \leq & 
\frac{v_n^2}{nh_n^d}\log 4+\max\left\{\frac{v_n^2}{nh_n^d}\log\mathbb{P}\left[\sqrt{v_n}\|B_n(x)\|>\delta\right] \ ; \ \frac{v_n^2}{nh_n^d}\log\mathbb{P}\left[\sqrt{v_n}\left|g(x)-g_n(x)\right| >\delta\frac{g(x)}{2}\right] \right\}.
\end{eqnarray*}
Now, since the sequence $\left(v_nB_n(x)\right)$ satisfies a LDP with speed $\left(\dfrac{nh_n^d}{v_n^2}\right)$ and good rate function $G_x$, there exists $c_1>0$ such that
\begin{eqnarray*}
\limsup_{n \to \infty}\frac{v_n}{nh_n^d}\log\mathbb{P}\left[\sqrt{v_n}\|B_n(x)\|>\delta\right]& \leq & -c_1.
\end{eqnarray*}
Moreover, the application of Theorem $1$ in Mokkadem et al. \cite{mokkapellwormsr} guarantees the existence of $c_2>0$ such that
\begin{eqnarray*}
\limsup_{n \to \infty}\frac{v_n}{nh_n^d}\log\mathbb{P}\left[\sqrt{v_n}\left|g(x)-g_n(x)\right| >\delta\frac{g(x)}{2}\right]<-c_2.
\end{eqnarray*}
We thus deduce that
\begin{eqnarray*}
\lim_{n \to \infty}\frac{v_n^2}{nh_n^d}\log\mathbb{P}\Big[v_n\|\left(r_n(x)-r(x)\right)-B_n(x)\|>\delta\Big]=-\infty,
\end{eqnarray*}
which means that the sequences $\left(v_n(r_n(x)-r(x))\right)$ and $\left(v_nB_n(x)\right)$ are exponentially contiguous. Theorem \ref{mdpregression} thus follows.  $\blacksquare$

\subsubsection{Proof of Lemma \ref{mdpbn}}
For any $ u \in \mathbb{R}^q$, set
\begin{eqnarray*}
\Gamma_{n,x}(u) & = & \frac{v_n^2}{nh_n^d}\log\mathbb{E}\left[\exp\left(\frac{nh_n^d}{v_n}\langle u,B_n(x)-\mathbb{E}\left(B_n(x)\right)\rangle\right)\right],\\
\Phi_x(u) & = & \frac{1}{2g^2(x)}\int_{\mathbb{R}^d\times\mathbb{R}^q}\langle u,y-r(x)\rangle^2 K^2(z)f(x,y)dzdy\\ & = & 
\frac{u^T\Sigma_x u}{2g(x)}\int_{\mathbb{R}^d}K^2(z)dz.
\end{eqnarray*}
To prove Lemma \ref{mdpbn}, it suffices to show that, for all $u \in \mathbb{R}^q$, 
\begin{eqnarray*}
\lim_{n \to \infty}\Gamma_{n,x}(u)& = & \Phi_x(u).
\end{eqnarray*}
As a matter of fact, since $\Phi_x$ is a quadratic function, Lemma \ref{mdpbn} then follows from the application of the G\"artner-Ellis Theorem.
For $u\in \mathbb{R}^q$, set 
\begin{eqnarray*}
\hat Z_i & = & \langle u,Y_i-r(x)\rangle K\left(\frac{x-X_i}{h_n}\right),
\end{eqnarray*}
and note that
\begin{eqnarray*}
\Gamma_{n,x}(u) & = & \frac{v_n^2}{nh_n^d}\log\mathbb{E}\left[\exp\left(\frac{1}{v_ng(x)}\sum_{i=1}^n\left[\hat Z_i-\mathbb{E}(\hat Z_i)\right]\right)\right].
\end{eqnarray*}
Since $(X_i,Y_i)$, $i=1,\dots,n$ are independent and identically distributed, it holds that
\begin{eqnarray*}
\Gamma_{n,x}(u) & = & \frac{v_n^2}{h_n^d}\log\mathbb{E}\left[e^{\frac{\hat Z_n}{v_ng(x)}}\right]-\frac{v_n}{h_n^dg(x)}\mathbb{E}(\hat Z_n).
\end{eqnarray*}
Now, we follow the same lines as in the proof of Lemma \ref{cum}.
A Taylor's expansion ensures that there exists $\hat c_n$ between $1$ and $\mathbb{E}\Big[e^{\frac{\hat Z_n}{v_ng(x)}}\Big]$ such that
\begin{eqnarray*}
\Gamma_{n,x}(u) 
& = & 
\frac{v_n^2}{h_n^d}\mathbb{E}\left[e^{\frac{\hat Z_n}{v_ng(x)}}-1-\frac{\hat Z_n}{v_ng(x)}\right]-\hat R_{n,x}^{(1)}(u)\\
& = &v_n^2\int_{\mathbb{R}^d\times\mathbb{R}^q}\left[e^{\frac{1}{v_ng(x)}\langle u,y-r(x)\rangle K\left(z\right)}-1-\frac{1}{v_ng(x)}\langle u,y-r(x)\rangle K\left(z\right)\right]f(x,y)dzdy \\ \mbox{} && -\hat R_{n,x}^{(1)}(u)+\hat R_{n,x}^{(2)}(u),
\end{eqnarray*}
with 
\begin{eqnarray*}
\lefteqn{\hat R_{n,x}^{(1)}(u)}\\ & = & \frac{v_n^2}{2\hat c_n^2h_n^d}\left(\mathbb{E}\left[e^{\frac{\hat Z_n}{v_ng(x)}}-1\right]\right)^2\\
\hat R_{n,x}^{(2)}(u) & = & v_n^2\int_{\mathbb{R}^d\times\mathbb{R}^q}\left[e^{\frac{1}{v_ng(x)}\langle u,y-r(x)\rangle K\left(z\right)}-1-\frac{1}{v_ng(x)}\langle u,y-r(x)\rangle K\left(z\right)\right]\Big[f(x-h_nz,y)-f(x,y)\Big]dzdy,
\end{eqnarray*}
and
\begin{eqnarray*}
\frac{1}{\hat c_n} & \leq & \max\left\{1,\frac{1}{\mathbb{E}\left(e^{\frac{\hat Z_n}{v_ng(x)}}\right)}\right\}. 
\end{eqnarray*}
Noting that 
\begin{eqnarray}
\label{alpha}
\left|\mathbb{E}\left(\frac{\hat Z_n}{v_ng(x)}\right)\right| 
 & \leq & \frac{h_n^d\|u\|}{v_ng(x)}\int_{\mathbb{R}^d}\left|K(z)\right|\left[\int_{\mathbb{R}^q}\|y\|f(x-h_nz,y)dy\right]dz \nonumber\\ \mbox{} & & +
\frac{h_n^d\|u\|\|r(x)\|}{v_ng(x)}\int_{\mathbb{R}^d}\left|K(z)\right|\left[\int_{\mathbb{R}^q}f(x-h_nz,y)dy\right]dz.
\end{eqnarray}
It follows from \eqref{borne} that, 
\begin{eqnarray*}
\mathbb{E}\left(\frac{\hat Z_n}{v_ng(x)}\right) \to 0.
\end{eqnarray*}
We deduce that there exists $c' \in \mathbb{R}^*_+$ such that 
\begin{eqnarray*}
\frac{1}{\hat c_n^2} \leq c',
\end{eqnarray*}
and thus, in view of \eqref{borne}, 
\begin{eqnarray*}
\hat R_{n,x}^{(1)}(u) 
& \leq & 
\frac{c'}{2}\frac{v_n^2}{h_n^d}\left(\int_{\mathbb{R}^d\times\mathbb{R}^q }\left[e^{\frac{1}{v_ng(x)}\langle u,y-r(x)\rangle K\left(\frac{x-s}{h_n}\right)}-1\right]f(s,y)dsdy\right)^2 \\ & \leq & \frac{c'}{2g^2(x)}h_n^d\left(\int_{\mathbb{R}^d\times\mathbb{R}^q}\left|\langle u,y-r(x)\rangle K(z)\right|e^{\left|\frac{1}{g(x)}\langle u,y-r(x)\rangle K(z)\right|}f(x-h_nz,y)dzdy\right)^2 \\ & \leq &
\frac{c'e^{\frac{2\|K\|_{\infty}}{g(x)}\|u\|\|r(x)\|}}{2g^2(x)}h_n^d\left(\int_{\mathbb{R}^d\times\mathbb{R}^q}\left|\langle u,y-r(x)\rangle K(z)\right|e^{\frac{\|K\|_{\infty}}{g(x)}\|u\|\|y\|}f(x-h_nz,y)dzdy\right)^2 \\ & \leq & Bh_n^d,
\end{eqnarray*}
where $B$ is a constant, so that
\begin{eqnarray*}
\lim_{n \to \infty}\left|\hat R_{n,x}^{(1)}(u)\right|=0.
\end{eqnarray*} 
On the other hand, since $\forall x \in\mathbb{R}$, $\displaystyle e^x-1-x  =  \frac{x^2}{2}+\frac{x^3}{6}d(x)$, with $\displaystyle d(x)\leq e^{|x|}$, 
we get
\begin{eqnarray}
\label{reste}
\lefteqn{\hat R_{n,x}^{(2)}(u)}\nonumber\\ & =& 
\frac{1}{2g^2(x)}\int_{\mathbb{R}^d\times\mathbb{R}^q}\langle u,y-r(x)\rangle^2 K^2(z)\left[f(x-h_nz,y)-f(x,y)\right]dzdy +\mathcal{R}_{n,x}(u),
\end{eqnarray}
with
\begin{eqnarray*}
\left|\mathcal{R}_{n,x}(u)\right| & \leq &  \frac{1}{6v_ng^3(x)}\int_{\mathbb{R}^d\times\mathbb{R}^q}\left|\langle u,y-r(x)\rangle^3 K^3(z)\right|e^{\frac{1}{g(x)}|\langle u,y-r(x)\rangle K(z)|}\left|f(x-h_nz,y)-f(x,y)\right|dzdy.
\end{eqnarray*}
It follows from \eqref{borne} that $\mathcal{R}_{n,x}$ converges to $0$. Applying then (A2) and (A3), we find 
\begin{eqnarray*}
\lim_{n\to \infty}\left|\hat R_{n,x}^{(2)}(u)\right|=0.
\end{eqnarray*}
Finally, we have 
\begin{eqnarray*}
\lefteqn{\Gamma_{n,x}(u)}\\ & = &  v_n^2\int_{\mathbb{R}^d\times\mathbb{R}^q}\left[e^{\frac{1}{v_ng(x)}\langle u,y-r(x)\rangle K\left(z\right)}-1-\frac{1}{v_ng(x)}\langle u,y-r(x)\rangle K\left(z\right)\right]f(x,y)dzdy-\hat R_{n,x}^{(1)}(u)\\ \mbox{} & &+\hat R_{n,x}^{(2)}(u)\\ & = & 
\Phi_x(u)-\hat R_{n,x}^{(1)}(u)+\hat R_{n,x}^{(2)}(u)+\hat R_{n,x}^{(3)}(u),
\end{eqnarray*}
with
\begin{eqnarray*}
\lefteqn{\hat R_{n,x}^{(3)}(u)}\\ & = & 
 v_n^2\int_{\mathbb{R}^d\times\mathbb{R}^q}\left[e^{\frac{1}{v_ng(x)}\langle u,y-r(x)\rangle K\left(z\right)}-1-\frac{\langle u,y-r(x)\rangle K\left(z\right)}{v_ng(x)}-\frac{\langle u,y-r(x)\rangle^2K^2\left(z\right)}{2v_n^2g^2(x)}\right]f(x,y)dzdy.
\end{eqnarray*}
By the majoration $\left|e^x-1-x-\frac{x^2}{2}\right|\leq \left|\frac{x^3}{6}d(x)\right|$, we get
\begin{eqnarray*}
\lefteqn{\left|\hat R_{n,x}^{(3)}(u)\right|}\\ & \leq & \frac{1}{6v_ng^3(x)}\int_{\mathbb{R}^d\times\mathbb{R}^q}\left|\langle u,y-r(x)\rangle^3K^3\left(z\right)\right|e^{\frac{1}{v_ng(x)}\left|\langle u,y-r(x)\rangle K\left(z\right)\right|}f(x,y)dzdy,
\end{eqnarray*}
and \eqref{borne} ensures that  
\begin{eqnarray*}
\lim_{n \to \infty}\left|\hat R_{n,x}^{(3)}(u)\right| & = & 0,
\end{eqnarray*}
which concludes the proof of Lemma \ref{mdpbn}. $\blacksquare$

\subsubsection{Proof of Lemma \ref{biaisregression}}
Observe that
\begin{eqnarray}
\label{egaliteregression}
\mathbb{E}\left(B_n(x)\right) & = & \frac{1}{g(x)}\left[\mathbb{E}\left(m_n(x)\right)-m(x)\right]-\frac{r(x)}{g(x)}\left[\mathbb{E}\left(g_n(x)\right)-g(x)\right].
\end{eqnarray}
Since
\begin{eqnarray*}
\mathbb{E}\left(m_n(x)\right)-m(x) & = &  \frac{1}{h_n^d}\mathbb{E}\left(Y_1K\left(\frac{x-X_1}{h_n}\right)\right) -m(x)\\ & = &  \frac{1}{h_n^d}\int_{\mathbb{R}^d\times\mathbb{R}^q}yK\left(\frac{x-z}{h_n}\right)f(z,y)dzdy-m(x)
 \\ & = & 
 \frac{1}{h_n^d}\int_{\mathbb{R}^d}m(z)K\left(\frac{x-z}{h_n}\right)dz-m(x)\\ & = & \int_{\mathbb{R}^d}K(y)\left[m(x-h_ny)-m(x)\right]dy,
\end{eqnarray*}
Assumptions (A5)i), (A5)iii) and a Taylor's expansion of $m$ of order $p$ ensure that
\begin{eqnarray}
\label{biais1regression}
\mathbb{E}\left(m_n(x)\right)-m(x) & = & O\left(h_n^p\right).
\end{eqnarray}
Similarly, we have 
\begin{eqnarray}
\label{biais2regression}
\mathbb{E}\left(g_n(x)\right)-g(x) & = & O\left(h_n^p\right).
\end{eqnarray}
We deduce from \eqref{egaliteregression}, \eqref{biais1regression}, and \eqref{biais2regression} that
\begin{eqnarray*}
\mathbb{E}\left(B_n(x)\right) & = & O\left(h_n^p\right),
\end{eqnarray*}
and thus Lemma \ref{biaisregression} follows from Assumption (A5)ii). $\blacksquare$

\subsection{Proof of Theorem \ref{mdprimregression}}
Set
\begin{eqnarray*}
\tilde B_n(x) & = & \frac{1}{g(x)}\left(\tilde m_n(x)-m(x)\right)-\frac{r(x)}{g(x)}\left(\tilde g_n(x)-g(x)\right),
\end{eqnarray*}
and, for any $u\in \mathbb{R}^q$, 
\begin{eqnarray*}
\tilde\Gamma_{n,x}(u) & = & \frac{v_n^2}{nh_n^d}\log\mathbb{E}\left[\exp\left(\frac{nh_n^d}{v_n}\langle u,\tilde B_n(x)-\mathbb{E}\left(\tilde B_n(x)\right)\rangle\right)\right],\\
\tilde \Phi_{a,x}(u) & = & \frac{1}{2(1+ad)g^2(x)}\int_{\mathbb{R}^d\times\mathbb{R}^q}\langle u,y-r(x)\rangle^2 K^2(z)f(x,y)dzdy\\ & = & 
\frac{1}{1+ad} \frac{u^T\Sigma_x u}{2g(x)}\int_{\mathbb{R}^d}K^2(z)dz.
\end{eqnarray*}
By following the steps of the proof of Lemma \ref{mdpbn} and by using the property \eqref{relregress}, we prove that
\begin{eqnarray}
\label{cumtierce}
\lim_{n \to \infty}\tilde\Gamma_{n,x}(u) &= & \tilde \Phi_{a,x}(u).
\end{eqnarray}
We first show how \eqref{cumtierce} implies Theorem \ref{mdprimregression}.
The function $\tilde \Phi_{a,x}$ being quadratic, the application of the G\"artner-Ellis Theorem then ensures that
\begin{eqnarray}
\label{mdpbnprim}
\mbox{the sequence} \left(v_n\left(\tilde B_n(x)-\mathbb{E}\left(\tilde B_n(x)\right)\right)\right)\mbox{ satisfies a LDP} \nonumber \\
\mbox{ with speed} \left(\frac{nh_n^d}{v_n^2}\right) \mbox{and good rate function}\  \tilde G_{a,x}.
\end{eqnarray}
Now, following the proof of Lemma \ref{biaisregression}, we have 
\begin{eqnarray*}
\mathbb{E}\left(\tilde m_n(x)\right)-m(x) & = &  
 \frac{1}{n}\sum_{i=1}^n\frac{1}{h_i^d}\int_{\mathbb{R}^d\times\mathbb{R}^q}yK\left(\frac{x-z}{h_i}\right)f(z,y)dzdy-m(x)\\ & = & \frac{1}{n}\sum_{i=1}^n\int_{\mathbb{R}^d}K(y)\left[m(x-h_iy)-m(x)\right]dy.
\end{eqnarray*}
Here again, Assumptions (A5)i), (A5)iii) and a Taylor's expansion of $m$ of order $p$ ensure that
\begin{eqnarray*}
\mathbb{E}\left(\tilde m_n(x)\right)-m(x) & = & O\left(\frac{1}{n}\sum_{i=1}^nh_i^p\right).\\
\end{eqnarray*}
and similarly, 
\begin{eqnarray*}
\mathbb{E}\left(\tilde g_n(x)\right)-g(x) & = & O\left(\frac{1}{n}\sum_{i=1}^nh_i^p\right),
\end{eqnarray*}
thus
\begin{eqnarray*}
v_n\mathbb{E}\left(\tilde B_n(x)\right) & =& O\left(\frac{v_n}{n}\sum_{i=1}^nh_i^p\right).
\end{eqnarray*}

\begin{itemize}
\item If $a p<1$, since $(h_n)$ varies regularly with exponent $(-a)$, we have, in view of \eqref{relregress} and Assumption (A5)ii),
\begin{eqnarray*}
\frac{v_n}{n}\sum_{i=1}^nh_i^p=O\left(\frac{v_n}{n}\left[nh_n^p\right]\right)=o(1).
\end{eqnarray*}
\item If $a p>1$, we have $\sum_ih_i^p<\infty$ and thus, since $v_n=o(nh_n^d)$, we get
\begin{eqnarray*}
\frac{v_n}{n}\sum_{i=1}^nh_i^p =  O\left(\frac{v_n}{n}\right)   =   o(1).
\end{eqnarray*}
\item In the case $a p=1$, let $\mathcal{L}$ be the slowly varying function such that $h_n=n^{-a}\mathcal{L}(n)$, and set $\varepsilon>0$ small enough. Since $a(p-\varepsilon) <1$, we have $h_n^p=o(h_n^{p-\varepsilon})$, and in view of \eqref{relregress} and (A4),
\begin{eqnarray*}
\dfrac{v_n}{n}\sum_{i=1}^nh_i^p & =&  o\left(v_nh_n^{p-\varepsilon}\right)  =  
o\left(nh_n^{d+p-\varepsilon}\right) \\  &= & o\left(n^{1-a(d+p-\varepsilon)}\left[\mathcal{L}(n)\right]^{d+p-\varepsilon}\right)\\ & = & o\left(n^{-a(d-\varepsilon)}\left[\mathcal{L}(n)\right]^{d+p-\varepsilon}\right)  =  o(1).
\end{eqnarray*}
\end{itemize}
We thus deduce that
\begin{eqnarray}
\label{esperanceprimregression}
\lim_{n \to \infty} v_n\mathbb{E}(\tilde B_n(x)) & =&0.
\end{eqnarray}

To conclude the proof of Theorem \ref{mdprimregression}, we follow the same lines as for the proof of Theorem \ref{mdpregression} (see Section 3.3), except that we apply \eqref{mdpbnprim} instead of Lemma \ref{mdpbn}, \eqref{esperanceprimregression} instead of Lemma  \ref{biaisregression}, and Theorem 1 in Mokkadem et al. \cite{mokkapellthiamr} instead of Theorem $1$ in Mokkadem et al. \cite{mokkapellwormsr}. $\blacksquare$\\
Let us now prove \eqref{cumtierce}.
For $u\in \mathbb{R}^q$, set 
\begin{eqnarray*}
T_i & = & \langle u,Y_i-r(x)\rangle K\left(\frac{x-X_i}{h_i}\right),
\end{eqnarray*}
and note that
\begin{eqnarray*}
\tilde\Gamma_{n,x}(u) & = & \frac{v_n^2}{nh_n^d}\log\mathbb{E}\left[\exp\left(\frac{h_n^d}{v_ng(x)}\sum_{i=1}^n\frac{1}{h_i^d}\left[T_i-\mathbb{E}(T_i)\right]\right)\right].
\end{eqnarray*}
Since $(X_i,Y_i)$, $i=1,\dots,n$ are independent and identically distributed, it holds that
\begin{eqnarray*}
\tilde\Gamma_{n,x}(u) & = & \frac{v_n^2}{h_n^d}\sum_{i=1}^n\log\mathbb{E}\left[e^{\frac{h_n^dT_i}{v_ng(x)h_i^d}}\right]-\frac{v_n}{ng(x)}\sum_{i=1}^n\frac{1}{h_i^d}\mathbb{E}(T_i).
\end{eqnarray*}
By Taylor expansion, there exists $c_{i,n}$ between $1$ and $\mathbb{E}\Big[e^{\frac{h_n^dT_i}{v_ng(x)}h_i^d}\Big]$ such that
\begin{eqnarray*}
\log\mathbb{E}\Big[e^{\frac{h_n^dT_i}{v_ng(x)}h_i^d}\Big] & = & \mathbb{E}\Big[e^{\frac{h_n^dT_i}{v_ng(x)}h_i^d}-1\Big]-\frac{1}{2c_{i,n}^2}\left(\mathbb{E}\Big[e^{\frac{h_n^dT_i}{v_ng(x)}h_i^d}-1\Big]\right)^2,
\end{eqnarray*}
and $\tilde\Gamma_{n,x}$ can be rewritten as 

\begin{eqnarray*}
\tilde\Gamma_{n,x}(u) 
& = & 
\frac{v_n^2}{nh_n^d}\sum_{i=1}^n\mathbb{E}\left[e^{\frac{T_ih_n^d}{v_ng(x)h_i^d}}-1\right]-\frac{v_n^2}{2nh_n^d}\sum_{i=1}^n\frac{1}{c_{i,n}^2}\left(\mathbb{E}\left[e^{\frac{T_ih_n^d}{v_ng(x)h_i^d}}-1\right]\right)^2-\frac{v_n}{ng(x)}\sum_{i=1}^n\frac{1}{h_i^d}\mathbb{E}(T_i).
\end{eqnarray*}
A Taylor expansion implies again that there exists $c'_{i,n}$ between $0$ and $\frac{T_ih_n^d}{v_ng(x)h_i^d}$  such that
\begin{eqnarray*}
\lefteqn{\mathbb{E}\left[e^{\frac{T_ih_n^d}{v_ng(x)h_i^d}}-1\right]}\\ & = & \frac{h_n^d}{v_ng(x)h_i^d}\mathbb{E}(T_i)+\frac{1}{2}\left(\frac{h_n^d}{v_ng(x)h_i^d}\right)^2\mathbb{E}(T_i^2)+\frac{1}{6}\left(\frac{h_n^d}{v_ng(x)h_i^d}\right)^3\mathbb{E}(e^{c'_{i,n}}T_i^3).
\end{eqnarray*}
Therefore,
\begin{eqnarray}
\label{xouli}
\lefteqn{\tilde\Gamma_{n,x}(u)} \nonumber\\ &= & 
\frac{1}{2g^2(x)}\frac{1}{nh_n^{-d}}\sum_{i=1}^n\frac{1}{h_i^d}\int_{\mathbb{R}^d\times\mathbb{R}^q}\langle u,y-r(x)\rangle^2 K^2(z)f(x, y)dzdy+\ddot{R}_{n,x}^{(1)}(u)+\ddot{R}_{n,x}^{(2)}(u),
\end{eqnarray}
with
\begin{eqnarray*}
\ddot{R}_{n,x}^{(1)}(u) & = & \frac{1}{6}\frac{h_n^{2d}}{v_ng^3(x)}\frac{1}{n}\sum_{i=1}^n\frac{1}{h_i^{2d}}\mathbb{E}(e^{c'_{i,n}}T_i^3)-\frac{v_n^2}{2nh_n^d}\sum_{i=1}^n\frac{1}{c_{i,n}^2}\left(\mathbb{E}\left[e^{\frac{T_ih_n^d}{v_ng(x)h_i^d}}-1\right]\right)^2, \\
\ddot{R}_{n,x}^{(2)}(u) &= & \frac{h_n^{d}}{2g^2(x)}\frac{1}{n}\sum_{i=1}^n\frac{1}{h_i^d}\int_{\mathbb{R}^d\times\mathbb{R}^q}\langle u,y-r(x)\rangle^2 K^2(z)\left[f(x-h_iz,y)-f(x, y)\right]dzdy.
\end{eqnarray*}
In view of \eqref{relregress}, the first term in the right-hand-side of \eqref{xouli} converges to $\tilde\Phi_{a,x}$.\\
It remains to prove that $\ddot{R}_{n,x}^{(1)}$ and $\ddot{R}_{n,x}^{(2)}$ converge to $0$. We have
\begin{eqnarray*}
\left|\mathbb{E}\left[\frac{h_n^dT_i}{v_ng(x)h_i^d}\right]\right| & \leq & \frac{h_n^d}{v_ng(x)h_i^d}\int_{\mathbb{R}^d\times\mathbb{R}^q}\left|\langle u,y-r(x)\rangle K\left(\frac{x-s}{h_i}\right)\right|f(s,y)dsdy \\ & \leq & 
\frac{h_n^d}{v_ng(x)}\int_{\mathbb{R}^d\times\mathbb{R}^q}\left|\langle u,y-r(x)\rangle K\left(z\right)\right|f(x-h_iz,y)dzdy.
\end{eqnarray*}
In view of \eqref{borne}, the integral is bounded, thus
\begin{eqnarray*}
\lim_{n \to \infty}\sup_{i\leq n}\mathbb{E}\left[\frac{h_n^dT_i}{v_ng(x)h_i^d}\right]  & = &  0,
\end{eqnarray*}
so that, there exists $c>0$ such that
\begin{eqnarray*}
\frac{1}{c_{i,n}^2} \leq c.
\end{eqnarray*}
Now, on the one hand, since $|e^t-1|\leq |t|e^{|t|}$, and in view of \eqref{borne} and \eqref{sup}, we have
\begin{eqnarray*}
\mathbb{E}\left|e^{\frac{h_n^dT_i}{v_ng(x)h_i^d}}-1\right| & \leq & 
\frac{h_n^d}{v_ng(x)}\int_{\mathbb{R}^d\times\mathbb{R}^q}\left|\langle u,y-r(x)\rangle K(z)\right|e^{c\left|\frac{1}{g(x)}\langle u,y-r(x)\rangle K(z)\right|}f(x-h_iz,y)dzdy \\ & \leq & B_1\frac{h_n^d}{v_ng(x)},
\end{eqnarray*}
where $B_1$ and $c$ are constants. We deduce that

\begin{eqnarray*}
\lim_{n \to \infty}\frac{v_n^2}{2nh_n^d}\sum_{i=1}^n\frac{1}{c_{i,n}^2}\left(\mathbb{E}\left[e^{\frac{T_ih_n^d}{v_ng(x)h_i^d}}-1\right]\right)^2 & = & 0.
\end{eqnarray*}
On the other hand,
\begin{eqnarray*}
\mathbb{E}\left[T_i^3e^{c'_{i,n}}\right] & \leq & \mathbb{E}\left[|T_i|^3e^{|c'_{i,n}|}\right] \\ & \leq & 
h_i^d\int_{\mathbb{R}^d\times\mathbb{R}^q}\left|\langle u,y-r(x)\rangle K(z)\right|^3e^{\frac{c}{g(x)}\left|\langle u,y-r(x)\rangle K(z)\right|}f(x-h_iz,y)dzdy\\ & \leq & 
B_2h_i^d,
\end{eqnarray*}
where $B_2$ is a constant. Thus,
\begin{eqnarray*}
\label{ban}
\left|\frac{h_n^{2d}}{6nv_ng^3(x)}\sum_{i=1}^n\frac{1}{h_i^{2d}}\mathbb{E}(e^{c'_{i,n}}T_i^3)\right|  & \leq & 
\frac{h_n^d}{6v_ng^3(x)}\frac{B_2}{nh_n^{-d}}\sum_{i=1}^nh_i^{-d}.
\end{eqnarray*}
Since $\lim_{n \to \infty}\dfrac{h_n^d}{v_n}=0$, \eqref{relregress} ensures that
\begin{eqnarray*}
\lim_{n \to \infty}\frac{h_n^{2d}}{6v_ng^3(x)}\frac{1}{n}\sum_{i=1}^n\frac{1}{h_i^{2d}}\mathbb{E}(e^{c'_{i,n}}T_i^3) & = & 0,
\end{eqnarray*}
 which proves that
\begin{eqnarray*}
\lim_{n \to \infty}\left|\ddot{R}_{n,x}^{(1)}(u)\right|=0.
\end{eqnarray*}
Finally, using \eqref{relregress}, (A2) and (A3), we have
\begin{eqnarray*}
\lefteqn{\lim_{n \to \infty} \ddot{R}_{n,x}^{(2)}(u)}\\ & = & \lim_{n \to \infty}\frac{1}{2g^2(x)}\frac{\sum_{i=1}^nh_i^{-d}}{nh_n^{-d}}\frac{1}{\sum_{i=1}^nh_i^{-d}}\sum_{i=1}^nh_i^{-d}\int_{\mathbb{R}^d\times\mathbb{R}^q}\langle u,y-r(x)\rangle^2 K^2(z)\left[f(x-h_iz,y)-f(x, y)\right]dzdy\\ & = & 0,
\end{eqnarray*}
which proves \eqref{cumtierce}. $\blacksquare$

\subsection{Proof of Propositions \ref{convex1} and \ref{convexprim}}
\subsubsection{Proof of Proposition \ref{convex1}}
\begin{itemize}
\item The strict convexity of $\Psi_x$ follows from its definition, since for any $\gamma\in ]0,1[$, and $\left(u,v\right)\neq \left(u',v'\right)$,
\begin{eqnarray*}
\Psi_x\left(\gamma\left(u,v\right)+(1-\gamma)\left(u',v'\right)\right) & = & \Psi_x\left(\left(\gamma u+(1-\gamma)u',\gamma v+(1-\gamma)v'\right)\right) 
 \\ & = & 
\int_{\mathbb{R}^d\times \mathbb{R}^q}\left(e^{\left[\langle \gamma u+(1-\gamma)u',y\rangle+\gamma v+(1-\gamma)v'\right]K(z)}-1\right)f(x,y)dzdy \\ & <&
\gamma\int_{\mathbb{R}^d\times \mathbb{R}^q}\left(e^{\left(\langle u,y\rangle+v\right)K(z)}-1\right)f(x,y)dzdy \\ \mbox{} && + (1-\gamma)\int_{\mathbb{R}^d\times \mathbb{R}^q}\left(e^{\left(\langle u',y\rangle+v'\right)K(z)}-1\right)f(x,y)dzdy,
\end{eqnarray*}
where the last inequality follows from the fact that $x\mapsto e^x$ is strictly convex. \\
Since $|e^t-1| \leq |t|e^{|t|}$ $\forall t \in \mathbb{R}$ and $K$ is bounded and integrable, \eqref{borne} imply that
\begin{eqnarray*}
\lefteqn{ \int_{\mathbb{R}^d\times \mathbb{R}^q}\left|\left(e^{\left(\langle u,y\rangle+v\right)K(z)}-1\right)f(x,y)\right|dzdy} \\ & \leq & 
\int_{\mathbb{R}^d\times \mathbb{R}^q}\left|\left(\langle u,y\rangle+v\right)K(z)\right|e^{\left|\left(\langle u,y\rangle+v\right)K(z)\right|}f(x,y)dzdy\\  & \leq & 
e^{\left|v\right|\|K\|_{\infty}}\|u\|\int_{\mathbb{R}^d}\left|K(z)\right|dz\int_{\mathbb{R}^q}\|y\|e^{\|u\|\|y\|\|K\|_{\infty}}f(x,y)dy \\ \mbox{} & & +e^{\left|v\right|\|K\|_{\infty}}|v|\int_{\mathbb{R}^d}\left|K(z)\right|dz\int_{\mathbb{R}^q}e^{\|u\|\|y\|\|K\|_{\infty}}f(x,y)dy<\infty,
\end{eqnarray*}
which ensures the existence of $\Psi_x$. \\
Next, set 
\begin{eqnarray*}
h_x(u,v,y,z)=\left[e^{\left(\langle u,y\rangle+v\right)K(z)}-1\right]f(x,y).
\end{eqnarray*}
Since $h_x$ is differentiable with respect to $(u,v)$ and
\begin{eqnarray*}
\nabla h_x(u,v,y,z) & = & \left(
\begin{array}{c}
ye^{\left(\langle u,y\rangle+v\right)K(z)}K(z)f(x,y)\\
e^{\left(\langle u,y\rangle+v\right)K(z)}K(z)f(x,y)
\end{array}
\right),
\end{eqnarray*}
using Assumption (A1) and \eqref{borne}, it can be seen that $\Psi_x$ is differentiable on $\mathbb{R}^q\times\mathbb{R}$. Since $\Psi_x$ is a smooth convex on $\mathbb{R}^q\times\mathbb{R}$, it follows that $\Psi_x$ is essentially smooth so that $I_x$ is a good rate function on $\mathbb{R}^q\times \mathbb{R}$ (see Dembo and Zeitouni \cite{dembzeitr}), which proves the first part of Proposition \ref{convex1}. \\Now, observe that $\overset{\circ}{\mathcal{D}}(\Psi_x)=\mathbb{R}^q\times\mathbb{R}$, and since $\Psi_x$ is strictly convex, it holds that the pair $\left(\overset{\circ}{\mathcal{D}}(\Psi_x),\Psi_x\right)$ is a convex function of Legendre type. It follows that $\left(\overset{\circ}{\mathcal{D}}(I_x),I_x\right)$ is a convex function of Legendre type (See Rockafellar \cite{rock}). Thus, Part 2 of Proposition \ref{convex1} follows from Theorem 26.5 of Rockafellar \cite{rock}.   
\item Let us now assume that $\lambda(S_-)=0$. Thus 
\begin{eqnarray*}
\Psi_x(u,v)= \int_{\mathbb{R}^d\times \mathbb{R}^q}\left(e^{u^TyK(z)}e^{vK(z)}-1\right)\mathds{1}_{S_+}(z)f(x,y)dzdy.
\end{eqnarray*}
For each $u\in\mathbb{R}^q$, the function $v\mapsto\left(e^{u^TyK(z)}e^{vK(z)}-1\right)\mathds{1}_{S_+}(z)f(x,y)$ is increasing in $v$ and goes to $-f(x,y)$ when $v\to -\infty$. Thus $\lim_{v \to -\infty}\Psi_x(u,v)=-g(x)\lambda (S_+)$ and $I_x(\vec{0},0)=g(x)\lambda (S_+)$.
Now, when $t_1\neq \vec{0}$, let us show that 
\begin{eqnarray*}
I_x(t_1,0)=+\infty.
\end{eqnarray*}
Let $M>0$, $\epsilon>0$ and set $u=(M+\epsilon)t_1/\|t_1\|^2$. Let $v\in \mathbb{R}$ such that
\begin{eqnarray*}
\left\lbrace
\begin{array}{ll}
-\Psi_x(u,v)>g(x)\lambda(S_+)-\epsilon \ \ \mbox{if}\ \ \lambda(S_+)<\infty\\
-\Psi_x(u,v)>M\ \ \mbox{if}\ \ \lambda(S_+)=\infty.
\end{array}
\right.
\end{eqnarray*}
Then, on the one hand, when $\lambda(S_+)<\infty$, we have
\begin{eqnarray*}
u^Tt_1-\Psi_x(u,v)\geq M+\epsilon+g(x)\lambda(S_+)-\epsilon>M.
\end{eqnarray*}
On the other hand, when $\lambda(S_+)=\infty$, we get
\begin{eqnarray*}
u^Tt_1-\Psi_x(u,v)\geq M+\epsilon+M>M.
\end{eqnarray*}
It follows that $\sup_{u,v}\left(u^Tt_1-\Psi_x(u,v)\right)=+\infty$. $\blacksquare$
\end{itemize}
\subsubsection{Proof of Proposition \ref{convexprim}}
Following the same lines of the proof of Proposition \ref{convex1}, we prove Proposition \ref{convexprim}. When $\lambda(S_-)=0$, for each $u\in\mathbb{R}^q$ and $s\in ]0,1]$, the map $v\mapsto s^{-ad}\left(e^{s^{ad}u^TyK(z)}e^{vK(z)}-1\right)\mathds{1}_{S_+}(z)f(x,y)$ is increasing in $v$ and goes to $-s^{-ad}f(x,y)$ when $v\to -\infty$. We deduce that $\lim_{v \to -\infty}\tilde\Psi_{a,x}(u,v)=-g(x)\lambda (S_+)\int_0^1s^{-ad}=-g(x)\lambda (S_+)/(1-ad)$ and $I_x(\vec{0},0)=g(x)\lambda (S_+)/(1-ad)$. $\blacksquare$ 
\subsection{Proof of Propositions \ref{barca} and \ref{barcaprim} }
\subsubsection{Proof of Proposition \ref{barca}}
\begin{description}
\item (i) Let us prove the first part of Proposition \ref{barca}.
\begin{itemize}
 \item If $\alpha <I_x(\vec{0},0)$, set 
\begin{eqnarray*}
 G=\{(a,b) \in \mathbb{R}^q\times \mathbb{R}, \ \ I_x(a,b)\leq \alpha\} \ \ \mbox{and} \ \  \hat{G}=\{(s,t) \in \mathbb{R}^q\times \mathbb{R},\ \ \hat{I}_x(s,t)\leq \alpha\}.
\end{eqnarray*}
We first show that $\hat{G}$ is a compact subset of $\mathbb{R}^q\times \mathbb{R}$.\\
First, observe that since $I_x$ is a good rate function, $G$ is a compact subset of $\mathbb{R}^q\times \mathbb{R}$. Let us define the following function
\begin{eqnarray*}
F:\mathbb{R}^q\times \mathbb{R} &\rightarrow & \mathbb{R}^q\times \mathbb{R}\\
(s,t) & \mapsto &  (st,t).
\end{eqnarray*}
Observe that $F$ is continuous and $\hat{G}=F^{-1}(G)$. We deduce that $\hat{G}$ is a closed subset of $\mathbb{R}^q\times \mathbb{R}$. \\
Now, let $(s_n,t_n)$ be a sequence of real numbers of $\hat G$, there exists $(x_n,y_n)\in G$ such that $(x_n,y_n)=F(s_n,t_n)=(s_nt_n,t_n)\in G$.\\ The compactness of $G$ on $\mathbb{R}^q\times \mathbb{R}$ ensures that there exists a sequence of real numbers  $(x_{n_k},y_{n_k})\in G$ such that $(x_{n_k},y_{n_k}) \to (x_0,y_0)$ as $k \to \infty$, where $(x_0,y_0) \in G$. Therefore, $(s_{n_k}t_{n_k},t_{n_k}) \to (x_0,y_0)$ as $k \to \infty$.\\
Noting that Condition (C) ensures that $\forall s \in \mathbb{R}^q$, $I_x(s,0)\geq I_x(\vec{0},0)> \alpha$ so that $(s,0) \notin G$. \\
It follows that
$y_0\neq 0$, and thus $t_{n_k} \to y_0$ and $s_{n_k}\to s_0$ as $k \to \infty$, where $s_0=x_0/y_0$. We deduce that $(s_{n_k},t_{n_k})\to(s_0,y_0)$ as $k \to \infty$, so that $(s_0,y_0)\in \hat{G}$. Thus $\hat{G}$ is a compact set.
Now we claim that the set $A=\left\{s,\ J(s)\leq \alpha\right\}$ is the image of $\hat{G}$ by the continuous map $\pi:(s,t)\mapsto s$, and thus it is a compact.\\
Indeed, clearly $\pi(\hat{G})\subset A$. For the opposite inclusion, consider $\alpha<\alpha'<I_x\left(\vec{0},0\right)$; the set $\hat{G}'=\left\{\hat{I}_x(s,t)\leq \alpha'\right\}$ is compact. Let $s_0\in A$, since $J(s_0)\leq \alpha$, we have $J(s_0)=\inf_{(s_0,t) \in\hat{G}'}\hat{I}_x(s_0,t)$; by compacity, there exists $t_0$ such that $J(s_0)=\hat{I}_x(s_0,t_0)$; $(s_0,t_0)\in \hat{G}$ and $\pi(s_0,t_0)=s_0$, thus $A\subset \pi(\hat{G})$.
\item If $\alpha \geq I_x(\vec{0},0)$, let $s \in \mathbb{R}^q$, we have
\begin{eqnarray*}
J(s) &\leq &I_x(st,t) \ \ \forall t\\
& \leq & I_x(\vec{0},0)\\
& \leq & \alpha.
\end{eqnarray*}
We deduce that $\mathbb{R}^q\subseteq \left\{ J(s)\leq \alpha\right\}$
and the second part of Proposition \ref{barca} (i) follows.
\end{itemize}
\item (ii) It is an obvious consequence of (i) and the definitions of $J$ and $J^*$.
\item (iii) Assume that $J^*(s)<\infty$.
\begin{itemize}
\item If $I_x(\vec{0},0)>\inf_tI_x(st,t)$, then 
\begin{eqnarray*}
\inf_tI_x(st,t)& = \inf_{t \neq 0}I_x(st,t),
\end{eqnarray*}
so that $J(s)=J^*(s)$.
\item If $I_x(\vec{0},0)=\inf_tI_x(st,t)$, since $J^*(s)<\infty$, there exists $t_0\neq 0$ such that $I_x(st_0,t_0)<\infty$. By the convexity of $I_x$, we have for any $\nu \in ]0,1], \ \ I_x(st_0\nu,t_0\nu) <\infty$ and
\begin{eqnarray*}
I_x(st_0\nu,t_0\nu) \leq \nu I_x(st_0,t_0)+(1-\nu)I_x(\vec{0},0). 
\end{eqnarray*}
We deduce that 
\begin{eqnarray*}
0\leq I_x(st_0\nu,t_0\nu)-I_x(\vec{0},0)\leq \nu\left(I_x(st_0,t_0)-I_x(\vec{0},0)\right),
\end{eqnarray*}
and if we take $\nu \to 0$, the third part of Proposition \ref{barca} follows.

\end{itemize}
\item (iv) Let us suppose that $\alpha<I_x(\vec{0},0)$ and let $s \in \left\{J^*(s)\leq \alpha \right\}$, then we have $J^*(s)<\infty$. We deduce from (iii) that $J(s)=J^*(s)$.
It follows that $J(s)\leq \alpha$, which ensures that $s \in \left\{ J(s)\leq \alpha \right\}$.\\
Conversely, if $s \in \left\{J(s)\leq \alpha \right\}$, then 
\begin{eqnarray*}
I_x(\vec{0},0)>\inf_{t}I_x(st,t),
\end{eqnarray*}
so that
\begin{eqnarray*}
\inf_{t}I_x(st,t) = \inf_{t\neq 0}I_x(st,t).
\end{eqnarray*}
That is $J(s)=J^*(s)$. Therefore, $J^*(s)\leq \alpha$, which ensures that
$s \in \left\{J^*(s)\leq \alpha \right\}$, and thus Proposition \ref{barca} is proved.
$\blacksquare$  
\end{description}
\subsubsection{Proof of Proposition \ref{barcaprim}}
Proposition \ref{barcaprim} is proved by following the same approach as for the proof of Proposition \ref{barca} with replacing $I_x$, $J$ and $J^*$ by $\tilde I_{a,x}$, $\tilde J_a$ and $\tilde J^*_a$ respectively. $\blacksquare$

\end{document}